\documentclass{amsart}

\usepackage{algorithm}
\usepackage{algpseudocode}

\usepackage{fullpage}

\usepackage{amsmath}
\usepackage{amsfonts}
\usepackage{amssymb}
\usepackage{amsthm}
\usepackage{mathtools} 
\usepackage{latexsym}

\usepackage{physics}

\usepackage{enumerate}
\usepackage{cancel}
\usepackage{cases}
\usepackage{empheq}
\usepackage{multicol}

\usepackage{graphics} 

\usepackage{wrapfig}

\usepackage{color}

\usepackage[backgroundcolor=gray!30,linecolor=black]{todonotes}
   
\usepackage{mathrsfs}
\usepackage{fontenc} 
\usepackage{inputenc} 

\usepackage{verbatim}

\theoremstyle{plain}

\theoremstyle{definition}

\theoremstyle{remark}

\numberwithin{equation}{section} 
\numberwithin{figure}{section}   

\usepackage[square,comma,numbers,sort&compress]{natbib}
\usepackage[colorlinks=true, pdfborder={ 00 0}]{hyperref}
\hypersetup{urlcolor=blue, citecolor=red}
\usepackage{url}

\usepackage{subcaption}

\newcommand{\nT}{\mathbb T}

\renewcommand{\abs}[1]{\left\lvert#1\right\rvert}

\renewcommand{\norm}[1]{\left\|#1\right\|}

\newcounter{my_counter}
\renewcommand{\div}{\nabla \cdot} 
\newcommand{\lap}{\Delta} 

\usepackage{placeins} 

\usepackage{cleveref}

\title{
An Assessment of Ensemble Kalman Filter and Azouani-Olson-Titi Algorithms for Data Assimilation: A Comparative Study}
\date{\today}

\author{Ning Ning}

\address{Department of Statistics, Texas A\&M University, College Station, TX, USA}
\author{Collin Victor}

\address{Department of Mathematics, Texas A\&M University, College Station, TX, USA}

\begin{document}

\begin{abstract}
Continuous data assimilation (CDA) is a method that continuously integrates observational data into a dynamical system to improve model accuracy in real-time. The AOT algorithm is one of the most widely used methods in CDA due to its efficiency in incorporating observational data to enhance model accuracy. However, no research to date has evaluated the performance of the AOT algorithm compared to the most widely used DA method, the ensemble Kalman filter (EnKF). Hence, in this paper, we conduct an extensive numerical examination to evaluate and compare these two algorithms for CDA problems with measurement error, addressing this gap. By analyzing the one-dimensional Kuramoto-Sivashinsky equation and the two-dimensional Navier-Stokes equation, which are central to many applications and representative in CDA problems, we found a significant computational advantage of the AOT algorithm.
\end{abstract}






\maketitle

\section{Introduction}

Data assimilation (DA) refers to a broad set of methodologies designed to enhance the precision of simulations by integrating observational data with physical models. The objective of data assimilation is to provide improved estimates of the state of a system, thereby enabling better predictions and assessments. 
Classical data assimilation techniques rely on the Kalman filter (KF), a foundational algorithm renowned in estimation and predictive analytics for its effectiveness in systems marked by dynamic fluctuations \citep{10.1115/1.3662552}. Its operation unfolds through a multi-stage process, involving prediction and update phases, wherein initial system states are forecasted using a pre-established model and subsequently refined with noisy and indirect measurements. The KF excels in accurately estimating both observable and unobservable parameters, proving indispensable in applications demanding dependable real-time information, such as robotics and navigation systems. However, it is computationally expensive for high-dimensional systems, assumes linear system dynamics and observation models, and relies heavily on accurate system and noise models.

The ensemble KF (EnKF) was first proposed by G. Evensen in 1994 \citep{evensen1994sequential} and later amended in 1998 \citep{houtekamer1998data}. It offers a versatile adaptation tailored for high-dimensional, nonlinear systems by employing a Monte Carlo implementation of the Bayesian update problem \citep{evensen1994sequential,houtekamer1998data}. Unlike its predecessor, the EnKF sidesteps the computationally intensive task of maintaining a covariance matrix by harnessing an ensemble of state vectors and leveraging the sample covariance computed from this ensemble. This creative adaptation incorporates the Bayesian framework to iteratively merge forecasts with incoming data, fortifying the filter’s predictive model with observable data over time. It plays a pivotal role in enhancing simulation accuracy in meteorological data by blending models with measurements to estimate the state of weather systems. It has revolutionized numerical weather and ocean predictions in atmospheric and oceanic modeling. The ability to manage large numbers of variables and to operate under the assumption of Gaussian distributed data has made it an influential component in ensemble forecasting and related fields. 

Traditional DA methods typically involve updating the model state at discrete intervals, processing batches of data to correct the model's forecast. In contrast, continuous DA (CDA) is a specific type of DA that updates the model state incrementally as new observations become available. This method combines the model's forecast with real-time data, continuously adjusting the model state to better match the observed reality. Rather than waiting for a complete forecast cycle to finish, CDA updates the model state as soon as new data is received. This continuous approach allows for real-time corrections and adjustments, making CDA particularly effective in systems where timely integration of data is crucial for maintaining accuracy, such as in meteorology and oceanography. CDA's ability to provide ongoing updates helps reduce the growth of errors and improve the reliability of forecasts in dynamic, chaotic systems.

CDA can be mathematically described by augmenting the original differential equations governing the system with additional terms that incorporate observational data. Consider a dynamical system described by the partial differential equation:
\begin{equation}
	u_t = F(u),
\end{equation}
where \(u(t,x)\) is the state variable, \(t\) represents time, \(x\) represents spatial coordinates, and \(F\) is a (potentially nonlinear and nonlocal) differential operator known from the model. The nudging method introduces an additional term to the right-hand side of the equation to nudge the model state towards the observations. This modified equation becomes:
\begin{equation}
\label{eqn:nudging}
 v_t = F(v) + \mu P(u - v),
\end{equation}
where  \(\mu\) is the nudging coefficient that is a positive parameter that controls the strength of the nudging term, $P$ is an operator that maps the model state to the observed state, and \(u\) is the observational data, representing the measurements obtained from the real system. For an overview of nudging methods, see, e.g., \cite{Lakshmivarahan_Lewis_2013} and the references therein.

The Azouani-Olson-Titi (AOT) algorithm \cite{Azouani_Olson_Titi_2014,Azouani_Titi_2014}, sometimes referred to as ``interpolated nudging'', has emerged as promising approaches to CDA, showcasing its potential in enhancing the precision of simulations. The AOT algorithm, also known as the CDA algorithm, is a recent data assimilation technique based on feedback control at the partial differential equation (PDE) level, aiming to provide improved estimates of a system's state, thereby enabling better predictions and assessments. This method integrates observational data with physical models, and it has shown promise in various applications such as downscaling the Bénard convection system through velocity measurements and analyzing continuous data assimilation schemes for the Navier-Stokes equations and stiff problems in science and engineering. 
The AOT algorithm modifies the original PDEs:
\begin{equation}
\label{eqn:AOT}
	 v_t = F(v) + \mu I_h(u - v).
\end{equation}
Here, \(I_h\) is a spatial interpolation operator with associated spatial length scale, $h$, that maps the model state \(v\) onto the observed data, extracting the components of the model state that correspond to the available observations and interpolating them to the entire domain. Here $h$ is given to be the spatial scale at which the observations are gathered. The use of interpolation in the AOT method is a significant departure from usual nudging methods, with crucial effects on implementation, convergence rates, and the amount of measurement data required.

The effectiveness of the AOT algorithm is underscored by its application in constructing increasingly accurate representations of physical processes, including fluid flow models, climate variability, and materials science. A large body of recent literature has built upon the AOT algorithm; see, e.g.,  \cite{Albanez_Nussenzveig_Lopes_Titi_2016,
Altaf_Titi_Knio_Zhao_Mc_Cabe_Hoteit_2015,
Bessaih_Olson_Titi_2015,
Biswas_Foias_Mondaini_Titi_2018downscaling,
Biswas_Hudson_Larios_Pei_2017,
Biswas_Martinez_2017,
Celik_Olson_Titi_2019,
Farhat_GlattHoltz_Martinez_McQuarrie_Whitehead_2019,
Farhat_Jolly_Titi_2015,
Farhat_Lunasin_Titi_2016abridged,
Farhat_Lunasin_Titi_2016benard,
Farhat_Lunasin_Titi_2016_Charney,
Farhat_Lunasin_Titi_2017_Horizontal,
Foias_Mondaini_Titi_2016,
Foias_Jolly_Kravchenko_Titi_2014,
GarciaArchilla_Novo_Titi_2018,
Gesho_Olson_Titi_2015,
GlattHoltz_Kukavica_Vicol_2014,
Ibdah_Mondaini_Titi_2018uniform,
Jolly_Martinez_Olson_Titi_2018_blurred_SQG,
Jolly_Martinez_Titi_2017,
Jolly_Sadigov_Titi_2015,
Larios_Pei_2017_KSE_DA_NL,
Lunasin_Titi_2015,
Markowich_Titi_Trabelsi_2016_Darcy,
Mondaini_Titi_2018_SIAM_NA,
Pei_2019,
Rebholz_Zerfas_2021_alg_nudge}.  
Computational experiments on the AOT algorithm and its variants have been conducted, including applications to the two-dimensional Navier-Stokes equation (2D NSE) \cite{Gesho_Olson_Titi_2015}, the 2D B\'enard convection equations \cite{Altaf_Titi_Knio_Zhao_Mc_Cabe_Hoteit_2015}, and the one-dimensional Kuramoto-Sivashinsky equation (1D KSE) \cite{Lunasin_Titi_2015,Larios_Pei_2017_KSE_DA_NL}. 
We note that while there has been much research on the AOT algorithm, we have yet to see any research directly comparing the EnKF and AOT, with the exception of \cite{price2021topics} (which has yet to appear in a peer-reviewed publication), and so the goal of this paper is to fill that gap.
We note that our work significantly differs from \cite{price2021topics} as we consider the standard nudging algorithm given in \cite{Azouani_Olson_Titi_2014} instead of the local version they consider, we consider different values of the parameters of EnKF (including the number of ensembles) which allows us to obtain convergence to much higher levels of precision, and we test both data assimilation algorithms on 1D and 2D equations. 

The remaining sections of the paper unfold as follows. 
In Section \ref{sec:algorithms}, we provide descriptions of the AOT algorithm for CDA with measurement error.
The stochastic EnKF for this setting is provided in Algorithm \ref{alg:EnKF}. In Section \ref{sec:1D}, we compare the performances of these two algorithms on the 1D KSE, thoroughly examining their properties under various circumstances. In Section \ref{sec:2D}, we compare their performances on the 2D NSE. These two models were chosen because they are central to many applications in engineering, meteorology, oceanography, and other fields involving fluid flow, and are representative in 1D and 2D studies of CDA. We focus on 1D and 2D models as they strike a balance between complexity and manageability, allowing researchers to capture essential dynamics without overwhelming computational demands, which is why they were chosen for comparative analysis of the EnKF and AOT algorithms. We conclude and provide final remarks in Section \ref{sec:Conclusion}. Throughout the paper, we use AOT and nudging interchangeably, but our focus remains consistently on AOT.

\section{Two Data Assimilation Algorithms}
\label{sec:algorithms}
The $F$ operator in the continuous dynamical system
\begin{align}
    u_t = F(u)
\end{align}
represents a potentially nonlocal and nonlinear differential operator, assumed to be known. While there are algorithmic variants that can recover aspects of the model if necessary, for our purposes here, we assume exact knowledge of the model. Specifically, $F$ will be provided such that for the 1D KSE tackled in Section \ref{sec:1D},
\begin{align}
F(u)=-u_{xxxx} - \lambda u_{xx} - uu_x
\end{align}
and for the 2D NSE handled in Section \ref{sec:2D},
\begin{align}F(u)=\nu \lap u - u\cdot \nabla u  - \nabla p + f,\end{align}
where the details are left to these sections. 
Note that there should be an initial condition for $u$, however the full state of this system is unknown, including the initial state. That is, we do not know the initial state of the system $u_0$; if we did, we could potentially evolve this system forward in time perfectly without using data assimilation. Instead, we have partial observations of the system, which we will use to attempt to recover the true solution.

The AOT algorithm in its general form is given by 
\begin{align}
    v_t = F(v) + \mu I_h(\tilde{u} - v),\\
    v(x,0) = v_0.\nonumber
\end{align}
Here, $I_h$ is a spatial interpolant with an associated length scale $h$.
Note that there is an initial condition $v_0$, which we specify. Typically we will utilize $v_0 = I_h(\tilde{u}_0)$ or $v_0 \equiv 0$ for simplicity. In this work, we take $I_h := P_M$, where $P_M$ as described in the previous section denoting the projection onto the lowest $M$ Fourier modes. We use $\tilde{u}$ instead of $u$ to indicate that the observations of $u$ may be contaminated with errors. 
The AOT algorithm for CDA with stochastically noisy data was proposed in \cite{Bessaih_Olson_Titi_2015}, using the following stochastic
evolution equation:
\begin{align}
    v_t = F(v) + \mu I_h(\tilde{u} - v)+\mu \xi,\\
    v(x,0) = v_0.\nonumber
\end{align}
Here, $\xi$ is a random error independent identically distributed of Gaussian type.

For our observations, we consider \begin{align}P_M\tilde{u}(t,x) = P_M\left(u(t,x) + \eta(t,x)\right),\end{align}
where $\eta \sim \mathcal{N}(0,\sigma_O)$, with $\sigma_O$ indicating the variance of the noise in the observations. Specifically, as we are considering Fourier mode projections for the observations, we will consider $\eta(t,x)$ to be given as a matrix with random complex coefficients, $\hat{\eta}_k$, indexed by the wavenumber, $k$, with $\Re \hat{\eta}, \Im \hat{\eta}\sim \mathcal{N}(0,\sigma_0)$ such that
\begin{align}\label{eq:eta}
    \eta(t,x) = \sum_{\abs{\vec{k}}=1}^M \hat{\eta}_{\vec{k}}(t) e^{i 2\pi \vec{k}\cdot x}.
\end{align}
As we are considering observations of $\mathbb{R}$-valued solutions, we stipulate that $\hat{\eta}_{-k} = \overline{\hat{\eta}_{k}}$, where $\overline{\cdot}$ denotes complex conjugation. We note that we will use the above formulation for the generation of both the noise in the observations and the manufactured noise used to inflate the forecast for EnKF. Without loss of generality we will write $\eta \sim N(0,\sigma)$ with the understanding that each component of the Fourier transform of $\eta$ is generated as Gaussian white noise with variance $\sigma$.

\begin{algorithm}[t!]
\caption{Stochastic EnKF (\cite{asch2016DA})}
\label{alg:EnKF}
\begin{algorithmic}[1] 
\Require Observational data $\tilde{u}_n:=P_M(u(t_n)) + \eta(t_n)$, where $\eta(t_n) \sim \mathcal{N}(0,\sigma_O)$,\\
Number of ensembles $K$ (chosen to be at least $m$ typically), \\
Observational operator $H:= P_M$ (projection onto the lowest $m$ Fourier modes), \\
Forward model $M$,\\
$\sigma_E$ (variance of ensemble generated noise),\\
$\sigma_I$ (variance of inflation noise).
\bigskip

\State Initialize the ensemble analysis, noise, and forecast terms.
\For{$k = 1$ \textbf{to} $K$}
    \State \vspace{-0.15cm}\[v_{1,k}^a \equiv 0, \quad u^{err}_{1,k} \sim \mathcal{N}(1,\sigma_{E}), \quad v_{1,k}^f = \tilde{u}_1 + u^{err}_{1,k}\]
\EndFor
\bigskip

\State Update the ensembles in the main time loop.
\For{$n = 1$ \textbf{to} $m$}
   \State Generate ensemble observations: for $k = 1 \ldots K$,
   $$u^{err}_{n,k} \sim \mathcal{N}(0,\sigma_E)$$

    \State Compute ensemble means: 
    $$\bar{v}_n^f = \frac{1}{K}\sum_{k=1}^{K} v^f_{n,k}, \quad \bar{u}^{err}_n = \frac{1}{K}\sum_{k=1}^K u^{err}_{n,k}, \quad\bar{v}^{H}_n = \frac{1}{K}\sum_{k=1}^K H(v_{n,k}^f)$$

    \State Compute normalized anomalies: for $k = 1 \ldots K$,
    $$[V^f]_{n,k} = \frac{v^f_{n,k} - \bar{v}^f_n}{\sqrt{K-1}},\quad [U^f]_{n,k} = \frac{H(v^f_{n,k}) - v^{err}_{n,k} - \bar{v}^H_n + \bar{u}^{err}_n}{\sqrt{K-1}}$$

    \State Compute the Kalman gain:
    $$\mathbf{K} = V^f(U^f)^T(U^f(U^f)^T)^{-1}$$

    \State Update the ensembles: for $k = 1 \ldots K$,
    $$v^a_{n,k} = v^f_{n,k} + \mathbf{K} (\tilde{u}_n - H(v^f_{n,k}))$$

    \State Update the forecast: for $k = 1 \ldots K$,
    $$v^f_{n,k} = M(v^a_{n,k})$$

    \State Inflate the forecast: for $k = 1 \ldots K$,
    $$v^f_{n,k} = v^f_{n,k} + \varepsilon, \quad\text{where }\varepsilon \sim \mathcal{N}(0,\sigma_I).$$
\EndFor
\end{algorithmic}
\end{algorithm}








The EnKF has proven to be highly efficient in a wide range of both academic and operational data assimilation problems with measurement error. It can be seen as a reduced-order Kalman Filter, which only handles error statistics up to the second order (i.e., means and covariances), which is why it is often referred to as a Gaussian filter. It is based on the concept of particles—a collection of state vectors called ensemble members. Instead of propagating large covariance matrices, the errors are represented by these scattered particles, with their variability representing the uncertainty of the system's state. The ensemble members are propagated by the model without any linearization, thereby avoiding the need to derive the tangent linear model and circumventing the approximations involved in linearization. Over the past decades years, it has become very popular in its many forms. For this study, we implemented one of its most widely used forms, the stochastic EnKF, as detailed in Algorithm \ref{alg:EnKF}.

In the implementation of both data assimilation schemes, we utilize the interpolation operator $I_h= P_M$, which projects onto the lowest $M$ Fourier modes. We decompose the error into three distinct components for better visualization: the observed error, the unobserved error, and the total error. The observed error is given by \begin{align}\norm{P_M(u - v)}_{L^2},\end{align}
representing the errors detectable in the observations, excluding the observational errors.
The unobserved error is given by
    \begin{align} \norm{Q(u-v)}_{L^2},\end{align}
where $Q = I - P_M$ denotes the projection orthogonal to $P_M$, projecting onto all Fourier modes larger than $M$. This error is designated as ``unobserved'' because it exists entirely in the length scales too small to be observed by our choice of $I_h$.
The total error is simply the sum of these quantities and is given simply by
    \begin{align} \norm{u - v}_{L^2}.\end{align}

\section{One-dimensional Kuramoto-Sivashinsky Equation}
\label{sec:1D}
In this section, we consider 1D KSE:
\begin{align}
    u_t + u_{xxxx} + \lambda u_{xx} + uu_x = 0,\\
    u(t,x)=u(t,x+2L) \quad \text{and} \quad u(0,x) = u_0,\nonumber
\end{align}
where $0 \leq x < 2L$ is the spatial coordinate, $t$ is time, $\lambda$ is a physical parameter, and $u_0$ is the initial condition. In the rest of the paper we will use $\lambda = \frac{1}{2}$.
Here, $u(t,x) \in \mathbb{R}$ represents the velocity. This PDE is known to exhibit chaotic behavior and serve as a 1D toy model for investigating the performance of various data assimilation algorithms. We focus on this regime as it is more computationally feasible to conduct tests with large numbers of ensembles. 

We equip the equation with periodic boundary conditions over the domain $[0, 32\pi)$, that is $u(t,x + 32\pi) = u(t,x)$. This specific domain size is chosen to ensure chaotic dynamics and a fully resolved Fourier spectrum. It's important to note that the turbulent effects in this equation are influenced by the size of the domain. The KSE was implemented numerically using a pseudospectral code with 2/3 dealiasing for the nonlinear term and an integrating factor method to handle the linear derivative terms. The equation was evolved forward using an explicit Euler scheme (see Algorithm \ref{alg:IF_Euler} for full details). 

We are simulating $u$ and $v$ with a uniform spatial mesh with spatial resolution $N = 2^8$. 
We note that as we are using spectral methods, we are calculating the spatial derivatives using Fourier transforms. In particular, the linear operator $L$ in Algorithm \ref{alg:IF_Euler} is computed on the Fourier wavemodes and used to calculate $\hat{u}_{n+1}$. That is, in this context, $L$ is given to be the fourth and second order derivative terms of $F$, which is given in the frequency domain by $L := -k^4 +\lambda k^2$, where $k$ is the index for the Fourier wavemodes. Here, $\hat{u}$ represents the coefficients for the Fourier transform, 
\begin{align}\hat{u}(t,k) := \int_{-\infty}^\infty u(t,x) e^{-2\pi i k\cdot x}dx. \end{align} 
For simplicity we suppress the wavemode index $k$ and index only by the time in our formulations in \Cref{alg:IF_Euler} with the understanding that the equations for $\hat{u}_{n+1}$ and $\hat{v}_{n+1}$ must be solved for each at each timestep for each integer wavemode $-N/2 < k \leq N/2$.  
To obtain the solution in physical space, $u(t,x)$, we simply take the inverse Fourier transform:
\begin{align} {u}(t,x)= \int_{-\infty}^\infty \hat{u}(t,k) e^{2\pi i k\cdot x}dk. \end{align}

\begin{algorithm}
    \caption{Integrating Factor Explicit Euler Scheme with Explicit Nudging Term}
\label{alg:IF_Euler}
\begin{algorithmic}[1] 
\Require Observational data $\tilde{u}_n:=P_M(u(t_n)) + \varepsilon(t_n)$, where $\varepsilon(t) \sim \mathcal{N}(0,\sigma_O)$,\\
Observational operator $P_M$ (projection onto the lowest $M$ Fourier modes), \\
Nonlinear function $N(u,t)$,\\
Linear function $L$ for integrating factor,\\
Timestep value $\Delta t$,\\
$u_n$ and $v_n$, the value of solutions $u$ and $v$ at previous iteration,\\
$\hat{u}_n = \mathcal{F}(u_n)$ and $\hat{v}_n = \mathcal{F}(v_n)$ represent the Fourier transforms of $u_n$ and $v_n$, respectively.
\bigskip

\State Forward time integration $\hat{u}_{n+1} = M(\hat{u}_n):$
\begin{align}\hat{u}_{n+1} =e^{L\Delta t}\hat{u}_n + \Delta t e^{L\Delta t} \mathcal{F}{N(u_n,t_n)}\end{align}

\State Nudging forward time integration
\begin{align}\hat{v}_{n+1} = {M(\hat{v}_n)} + \Delta t e^{L\Delta t} \mu  P_M\left( \hat{\tilde{u}}_n - \hat{v}_n\right)
\end{align}

\end{algorithmic}
\end{algorithm}

\subsection{Computational Results}
\label{sec:1d_results}
\begin{figure}
    \centering
    \includegraphics[width=0.75\linewidth]{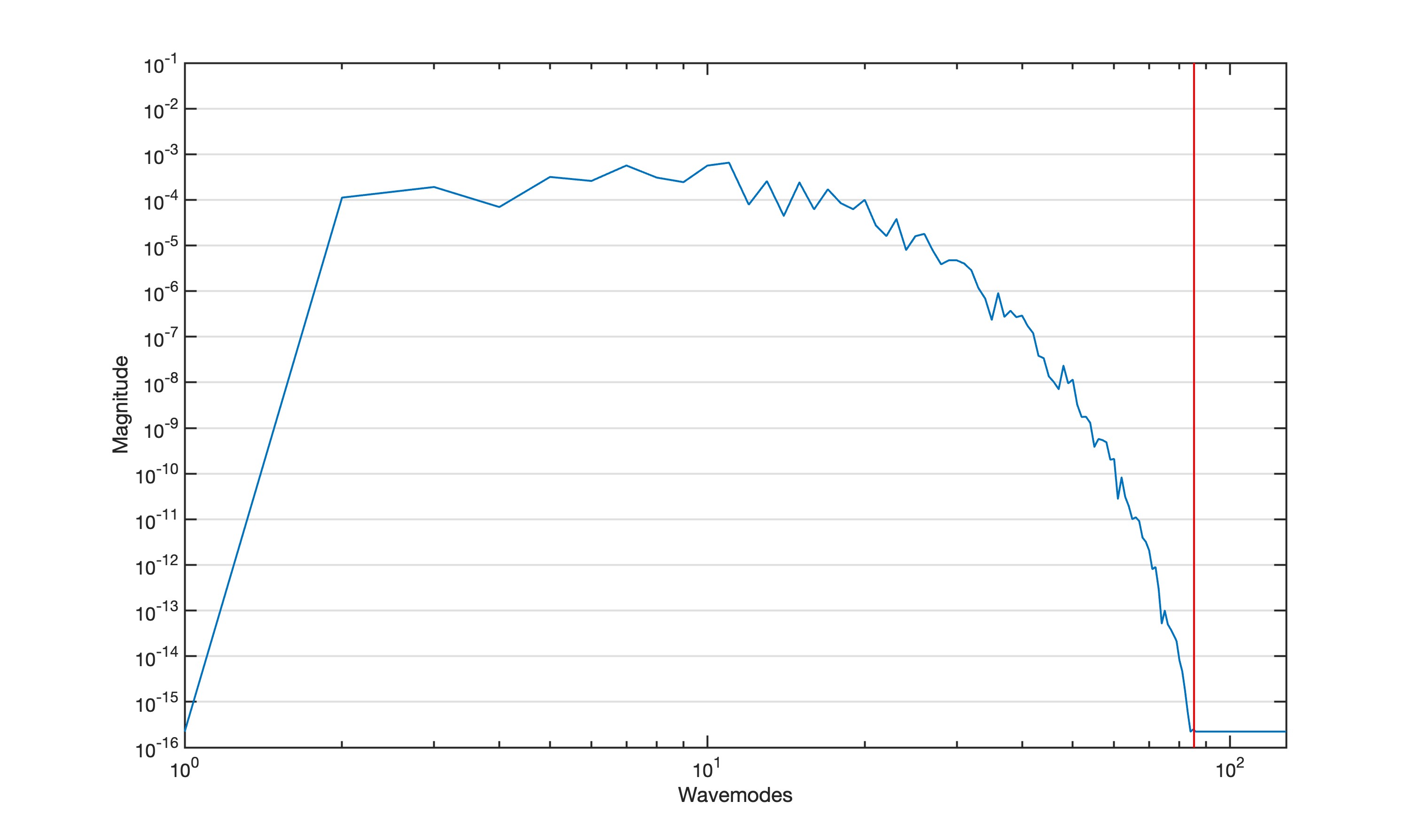}
    \caption{Energy spectrum of the initial data with $\lambda = 0.5$,  $\Delta t = 0.01$, on the spatial domain $[0, 32\pi)$. The vertical red line is the 2/3 dealiasing cutoff as $\frac{2}{3}\frac{N}{2} = 85.\overline{3}.$}
    \label{fig:spectrum KSE}
\end{figure}

To examine these data assimilation algorithms, we implemented an ``identical twin" experimental design. Such experiments are standard in the testing of data assimilation algorithms and involve running two separate simulations. The first simulation produces a reference solution, which we use to generate observational data and as the ``true" solution we attempt to recover using data assimilation. Next, we start the twin simulation, a second simulation that begins from a different set of initial conditions. The observational data is then assimilated into the twin simulation, and the error between both simulations is calculated. The phrase ``identical twin" refers to the underlying physics in the simulations; namely, both the reference and twin simulations utilize the same underlying model with consistent physics between them. For a detailed discussion of twin experiments, see, e.g., \cite{yu2019twin}.
In the current setting, we initialized a reference solution $u(t,x)$ using the initial data $u_0 = \cos(x/16)(1+\sin(x/16))$. We evolved this initial data up to time $t = 10,000$ to establish chaotic dynamics. We can see in \Cref{fig:spectrum KSE} that the initial profile is well-resolved, as the energy spectrum of the wavemodes decays to machine precision (approximately 2.2204e-16) by the dealiasing cutoff line. We note that the energy spectrum of all of the simulations featured in this section remain well-resolved for the duration of each simulation. We employed data assimilation methods to endeavor to recover this solution using only partial observations, with the aid of this reference solution. The partial observations in this setting are given by the first $16$ Fourier modes.

\begin{figure}[t!]
    \begin{subfigure}[b]{.32\textwidth}
\centering

            \includegraphics[width=\textwidth]{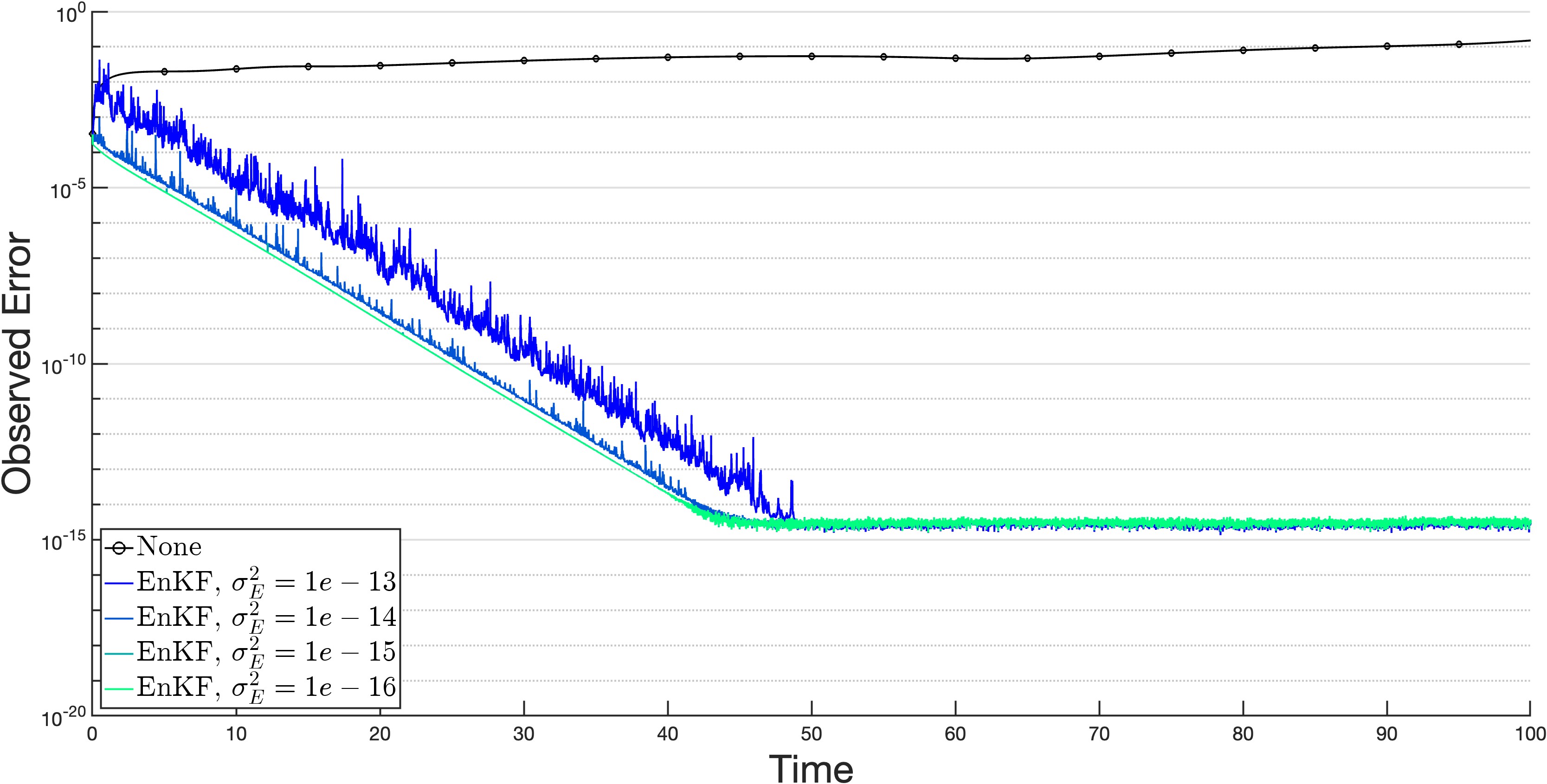}

    \end{subfigure}
    \begin{subfigure}[b]{.32\textwidth}
\centering

            \includegraphics[width=\textwidth]{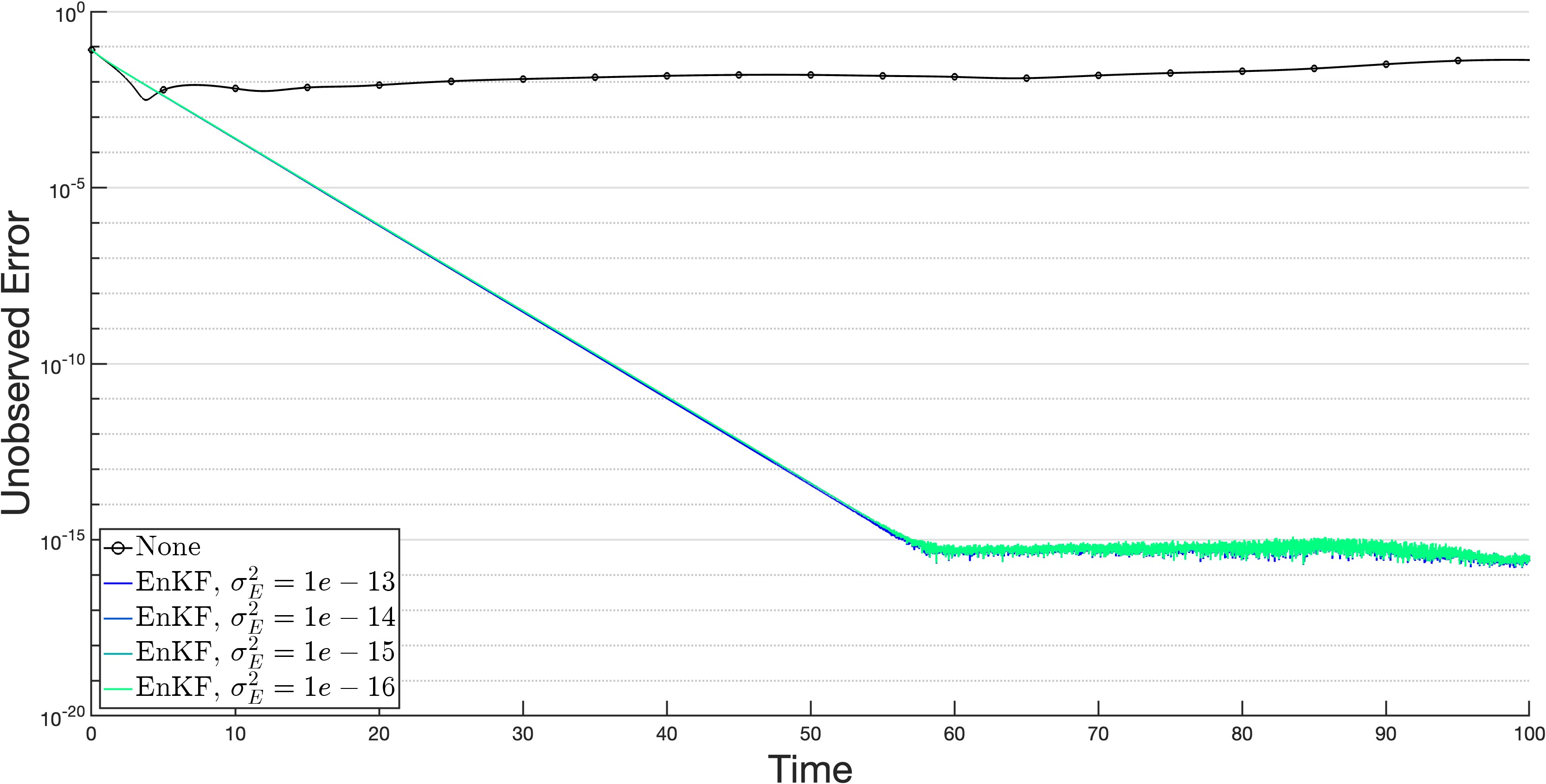}

    \end{subfigure}
    \begin{subfigure}[b]{.32\textwidth}
\centering

            \includegraphics[width=\textwidth]{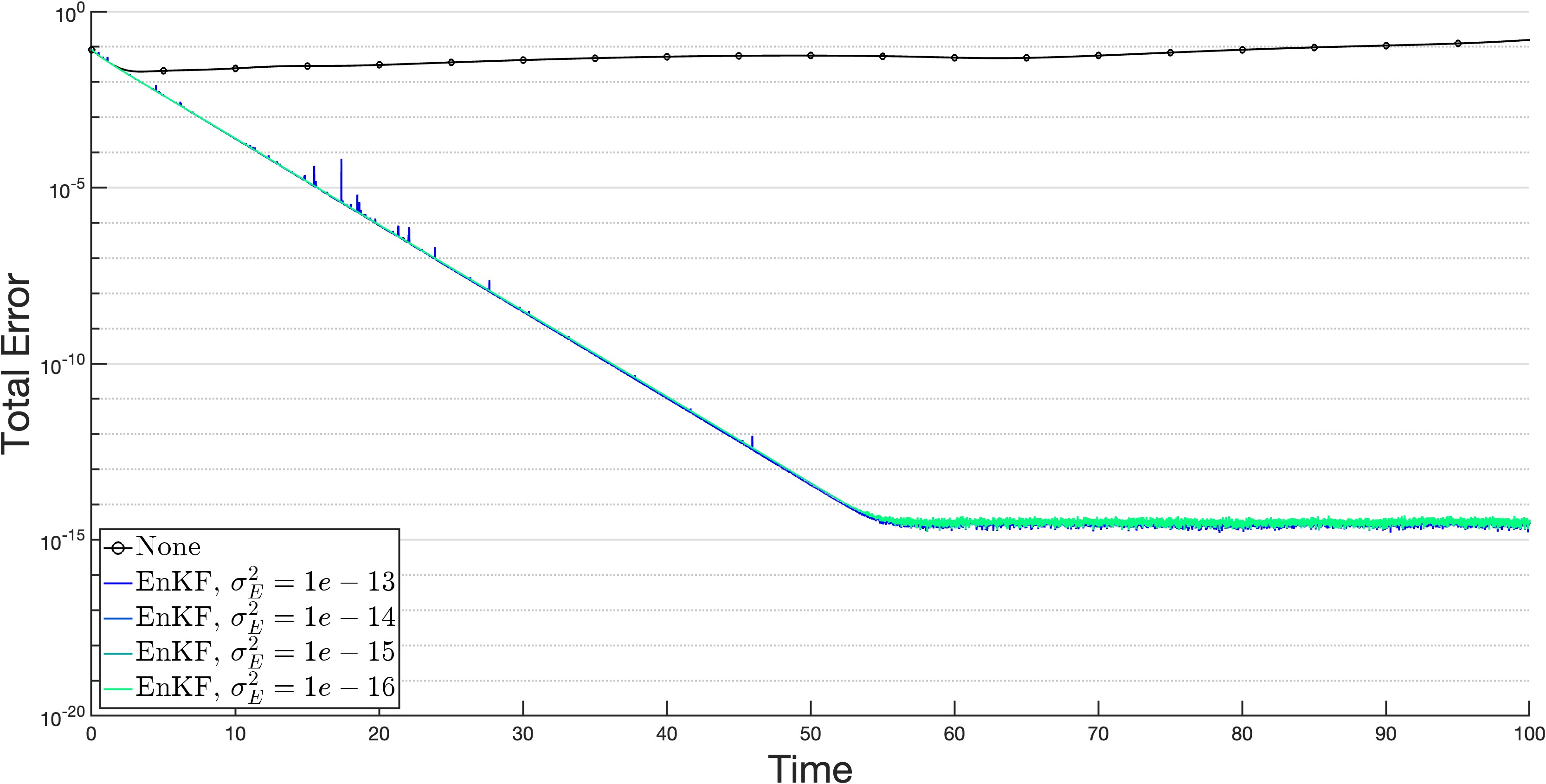}

    \end{subfigure}
    \caption{Error over time for different values of $\sigma_E$. Simulations utilized $K = 500$ ensemble members, $\sigma_O^2 = 0$, and $\sigma_I^2 = 1e-14$.}\label{fig:Ensemble_variance}
\end{figure}

In our simulations, we observed that adjusting the value of $\sigma_E$ had minimal impact on convergence for small values ($\sigma_E^2 \leq 1e-14$). As illustrated in Figure \ref{fig:Ensemble_variance} decreasing $\sigma_E$ seemed to reduce the error on the observed modes without significantly affecting the overall error on the rest of the system. It's worth mentioning that Figure \ref{fig:Ensemble_variance} depicts results with $K= 500$ ensemble members, while results for larger $\sigma_E$ values are not shown here. We found that the simulations were highly sensitive to the number of ensemble members, with larger values of $K$ leading to greater stability for larger $\sigma_E$ values. Given the decreased error and improved stability observed for small values of $\sigma_E$, we opted to use  $\sigma_E^2 = 1e-16$, which is approximately machine precision. We note that setting $\sigma_E$ in this manner poses no issue, as it is a user-controlled parameter.

\begin{figure}
    \begin{subfigure}[b]{.32\textwidth}
\centering

            \includegraphics[width=\textwidth]{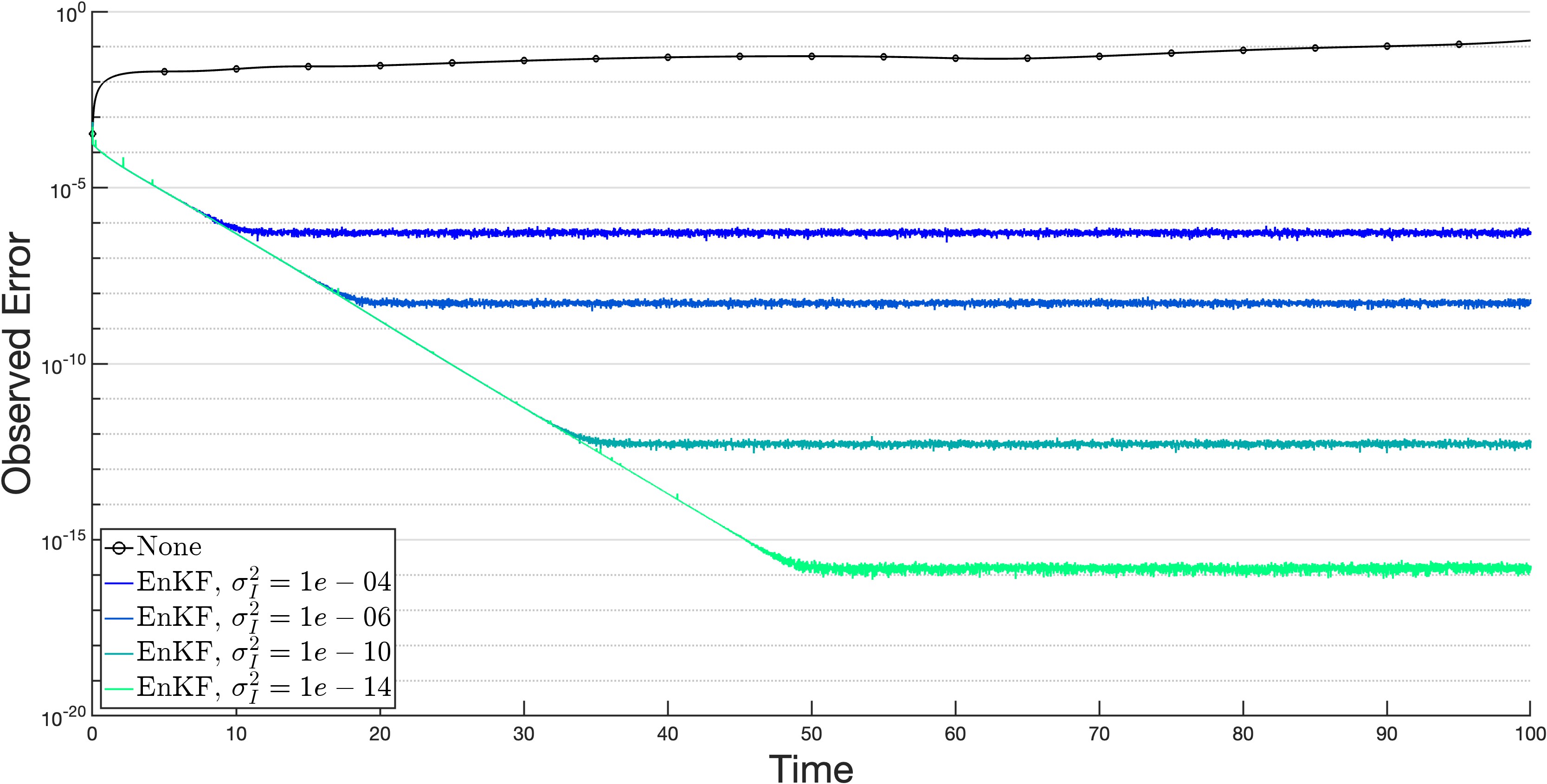}

    \end{subfigure}
    \begin{subfigure}[b]{.32\textwidth}
\centering

            \includegraphics[width=\textwidth]{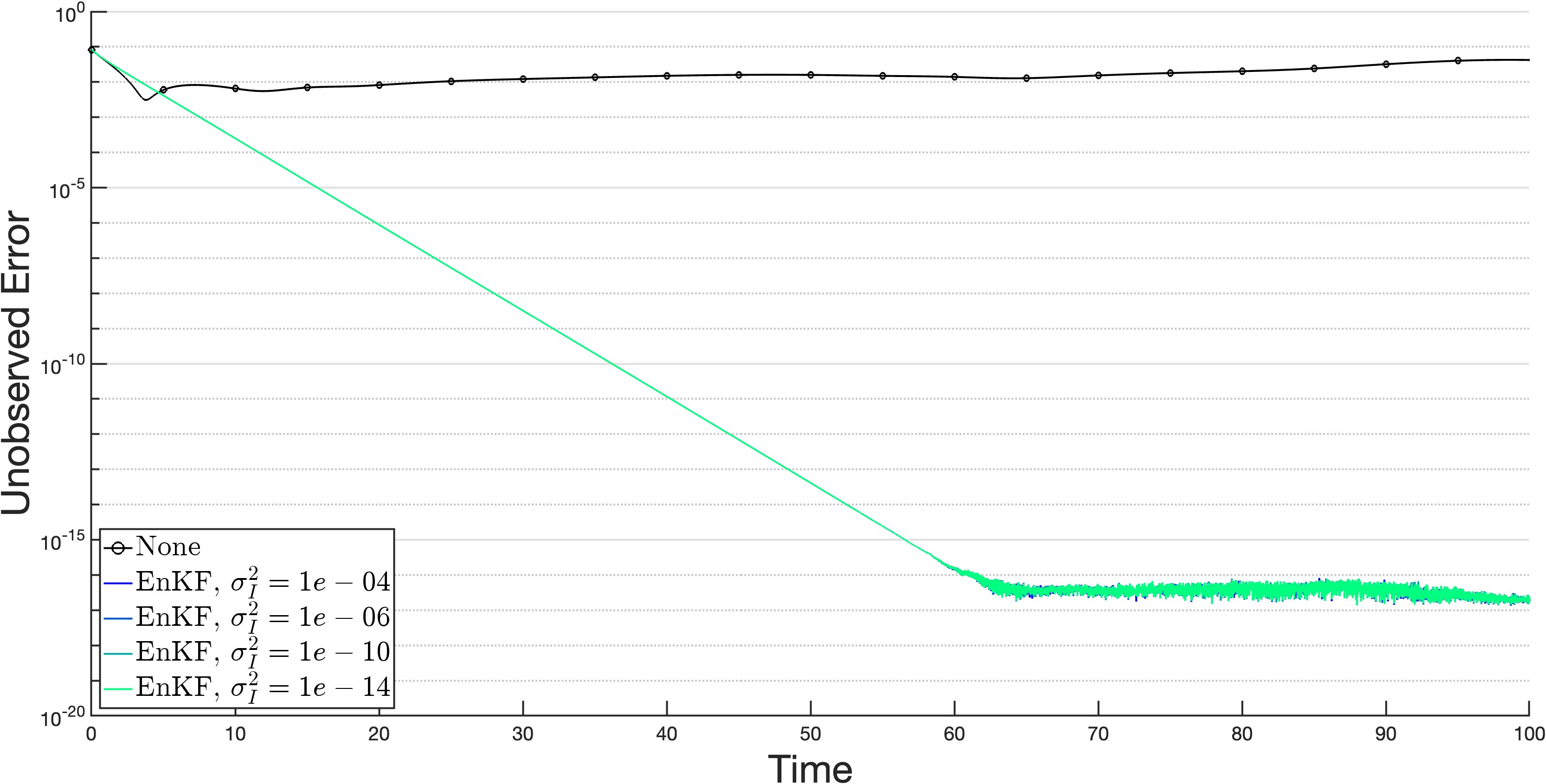}

    \end{subfigure}
    \begin{subfigure}[b]{.32\textwidth}
\centering

            \includegraphics[width=\textwidth]{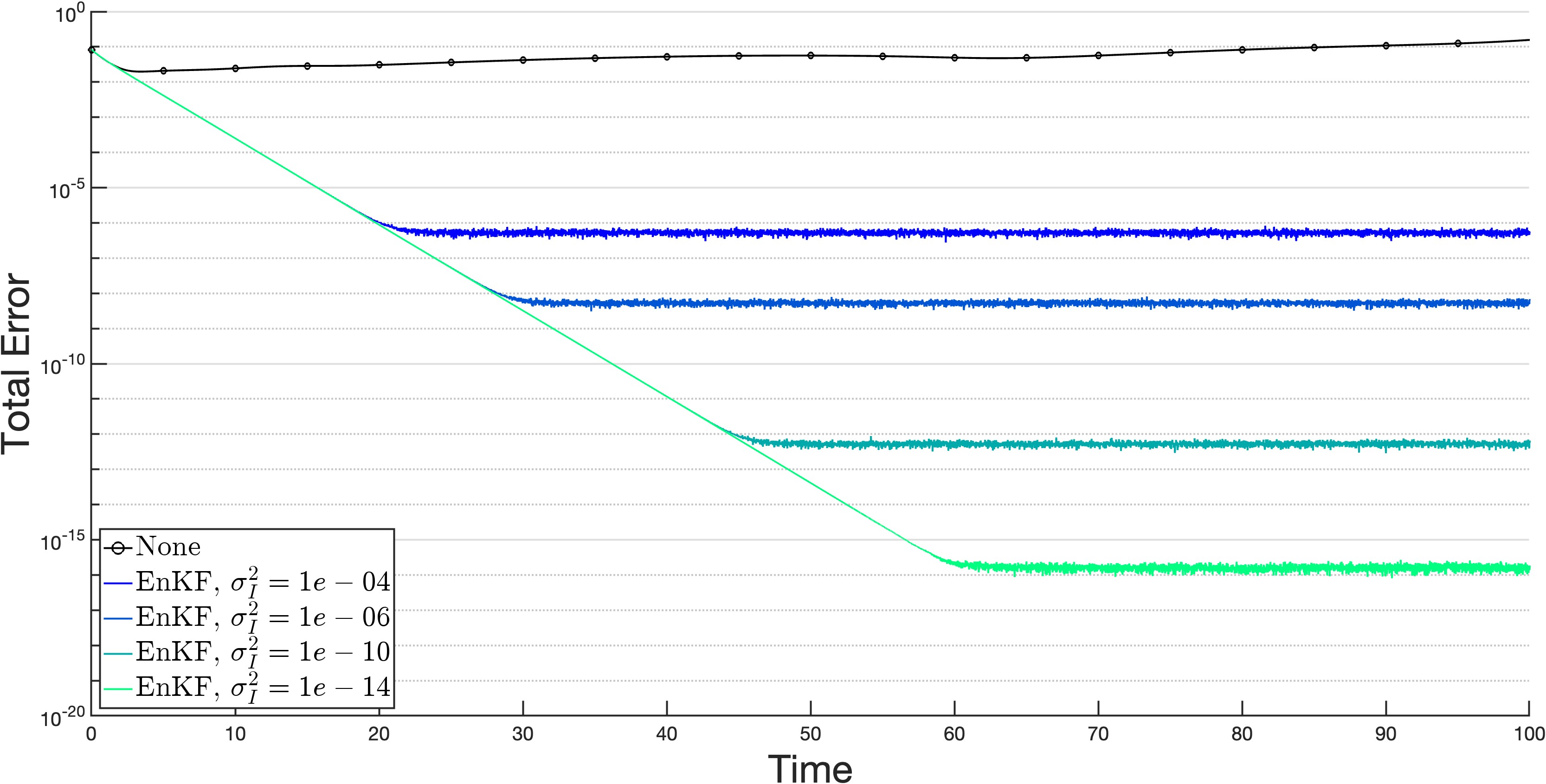}

    \end{subfigure}
    \caption{Error over time for different values of $\sigma_I$. Simulations utilized $K = 32$ ensemble members, $\sigma_O^2 = 0$, and $\sigma_E^2 = 1e-16$.}\label{fig:Covariance_Inflation}
\end{figure}

We observed that adjusting the value of $\sigma_I$ yielded intriguing results concerning the overall error, as depicted in Figure \ref{fig:Covariance_Inflation}. In the case of deterministic observations, we found that the error in the observed low modes converged to a fixed precision determined by $\sigma_I$, with smaller values leading to higher precision levels. Surprisingly, the lack of convergence seemed to have no discernible effect on the unobserved high modes, as convergence to machine precision in the ensemble average was consistently achieved. This observation appears counterintuitive, as one might expect errors in the low modes to propagate to the high modes. However, it seems that an averaging effect with the inflationary noise may be at play. Notably, the error on the low modes is calculated before the covariance inflation term is added to all ensemble members, indicating that this discrepancy isn't an implementation issue, as one might initially assume. To obtain the ``true" error, one should discount the error in the low modes as calculated using the low modes of the ensembles, and instead directly utilize the error between the observations and the truth, characterized by $\sigma_O$, along with the error in the unobserved modes of the ensemble solution. It's important to highlight that this issue of covariance inflation is absent in the AOT algorithm.

\begin{figure}
    \begin{subfigure}[b]{.32\textwidth}
\centering

            \includegraphics[width=\textwidth]{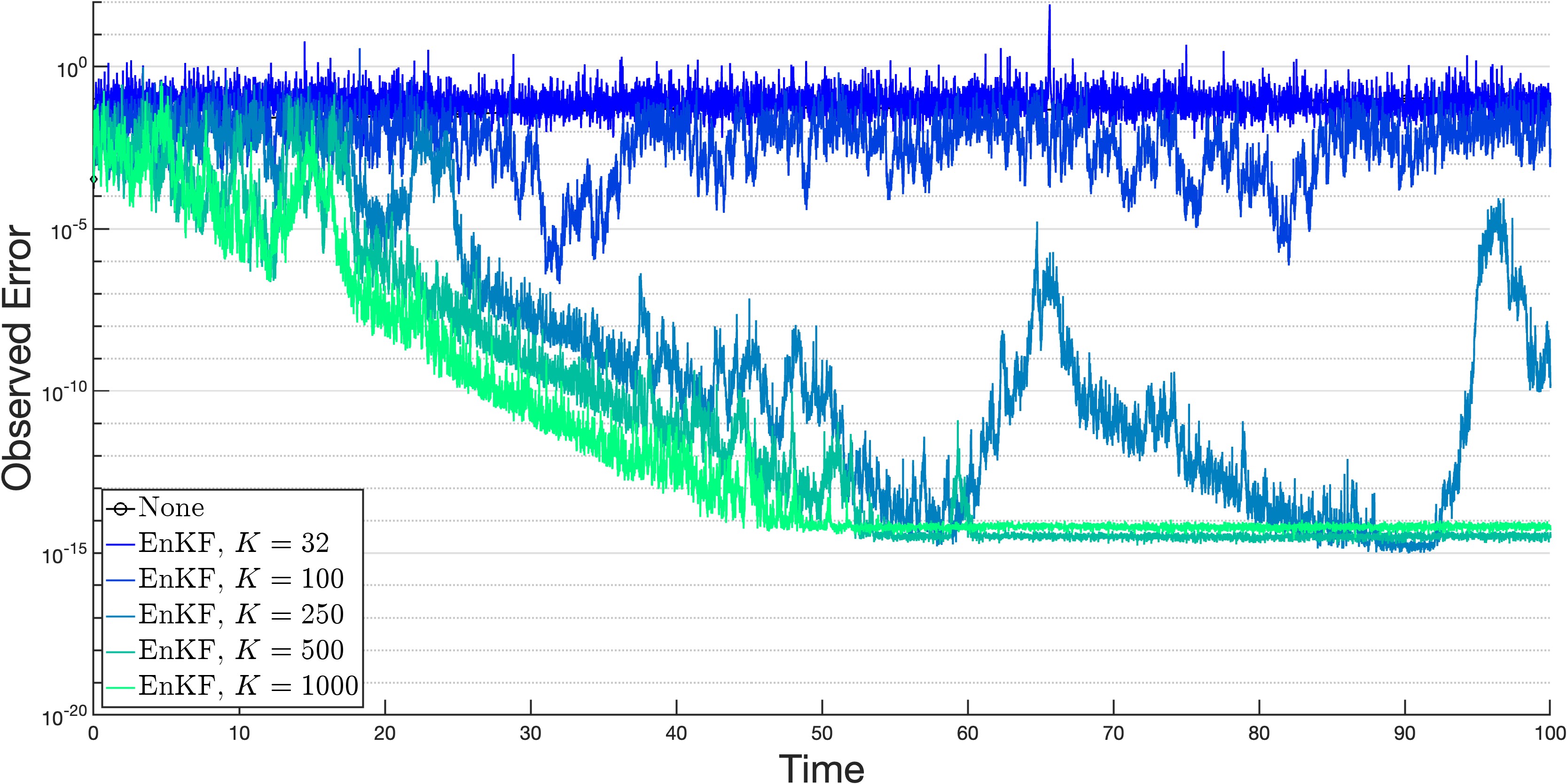}

    \end{subfigure}
    \begin{subfigure}[b]{.32\textwidth}
\centering

            \includegraphics[width=\textwidth]{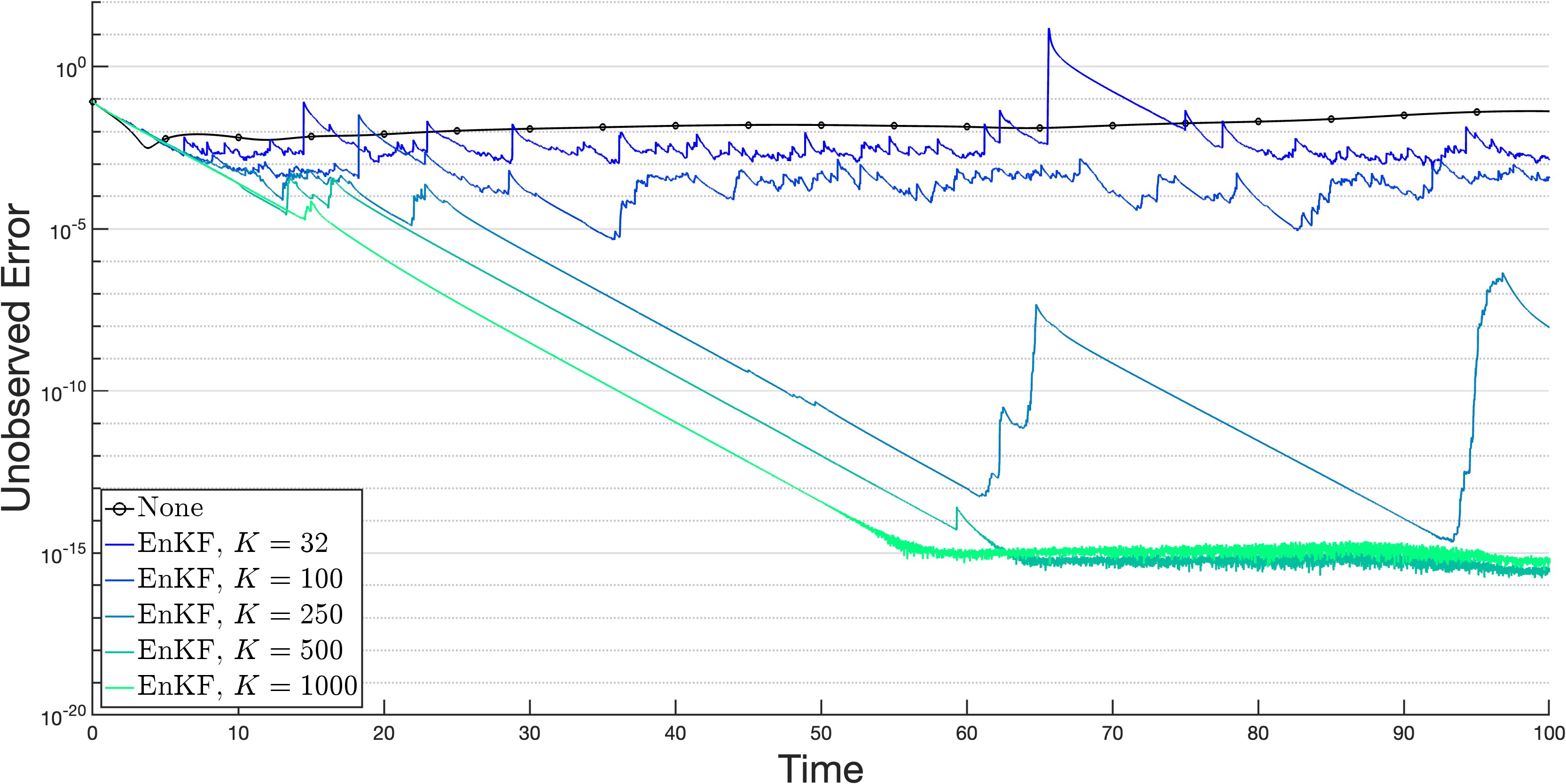}

    \end{subfigure}
    \begin{subfigure}[b]{.32\textwidth}
\centering

            \includegraphics[width=\textwidth]{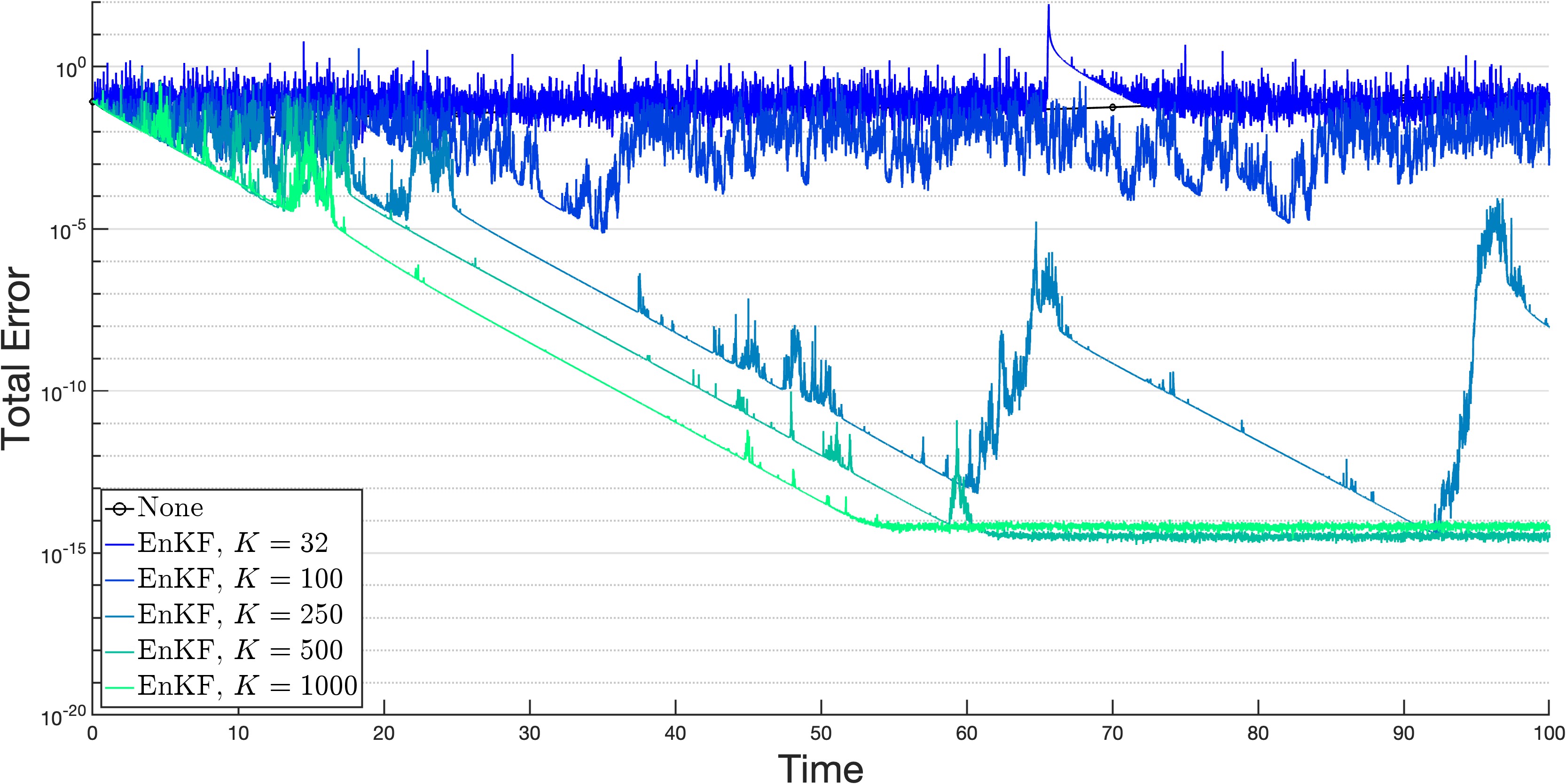}

    \end{subfigure}
    \caption{Error over time for different values of $K$ without covariance inflation. Simulations utilized $\sigma_I^2 = 0$, $\sigma_O^2 = 0$, and $\sigma_E^2 = 1e-16$.}\label{fig:Ensembles_no_inflation}
\end{figure}

\begin{figure}
    \begin{subfigure}[b]{.32\textwidth}
\centering

            \includegraphics[width=\textwidth]{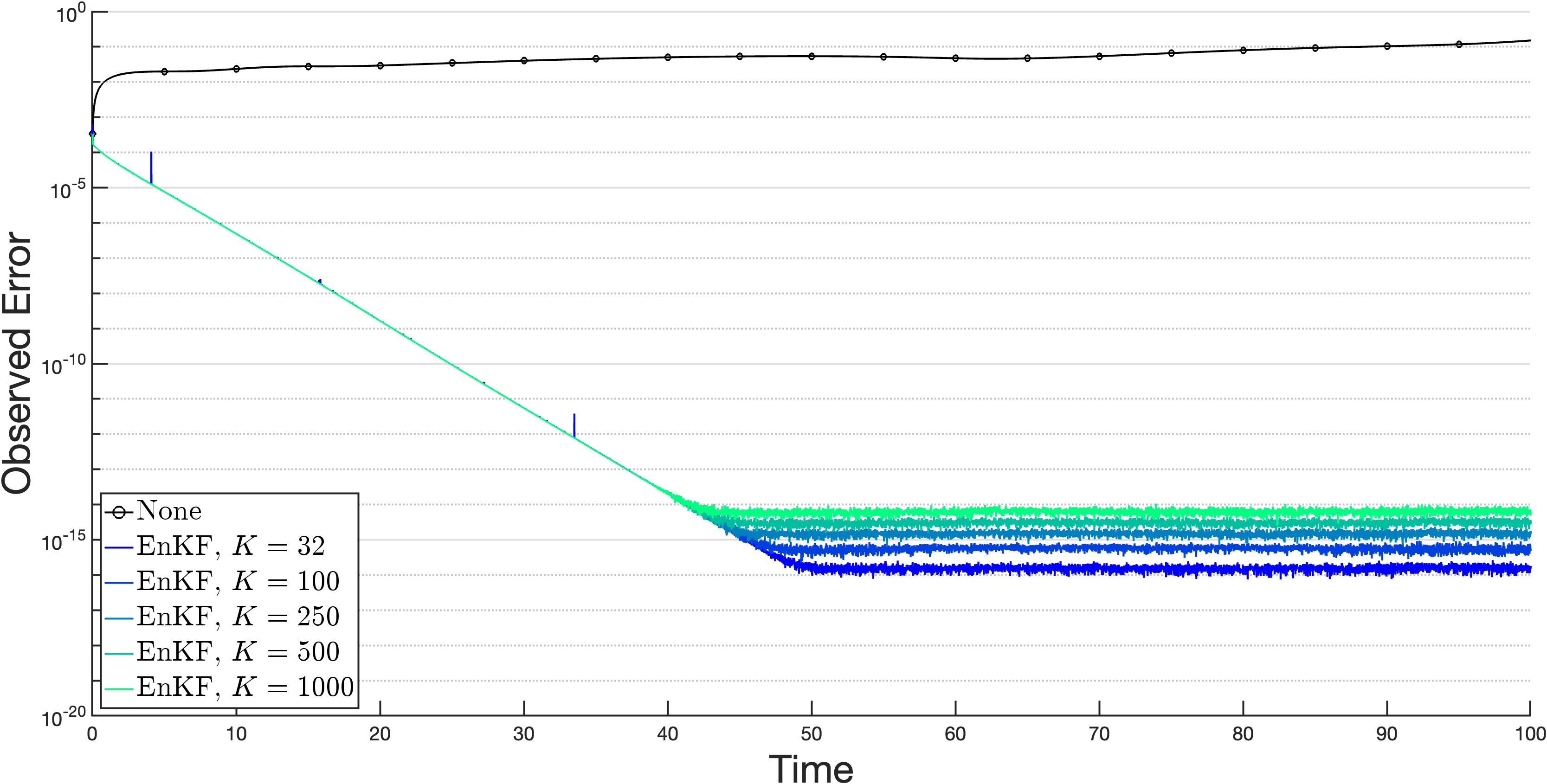}

    \end{subfigure}
    \begin{subfigure}[b]{.32\textwidth}
\centering

            \includegraphics[width=\textwidth]{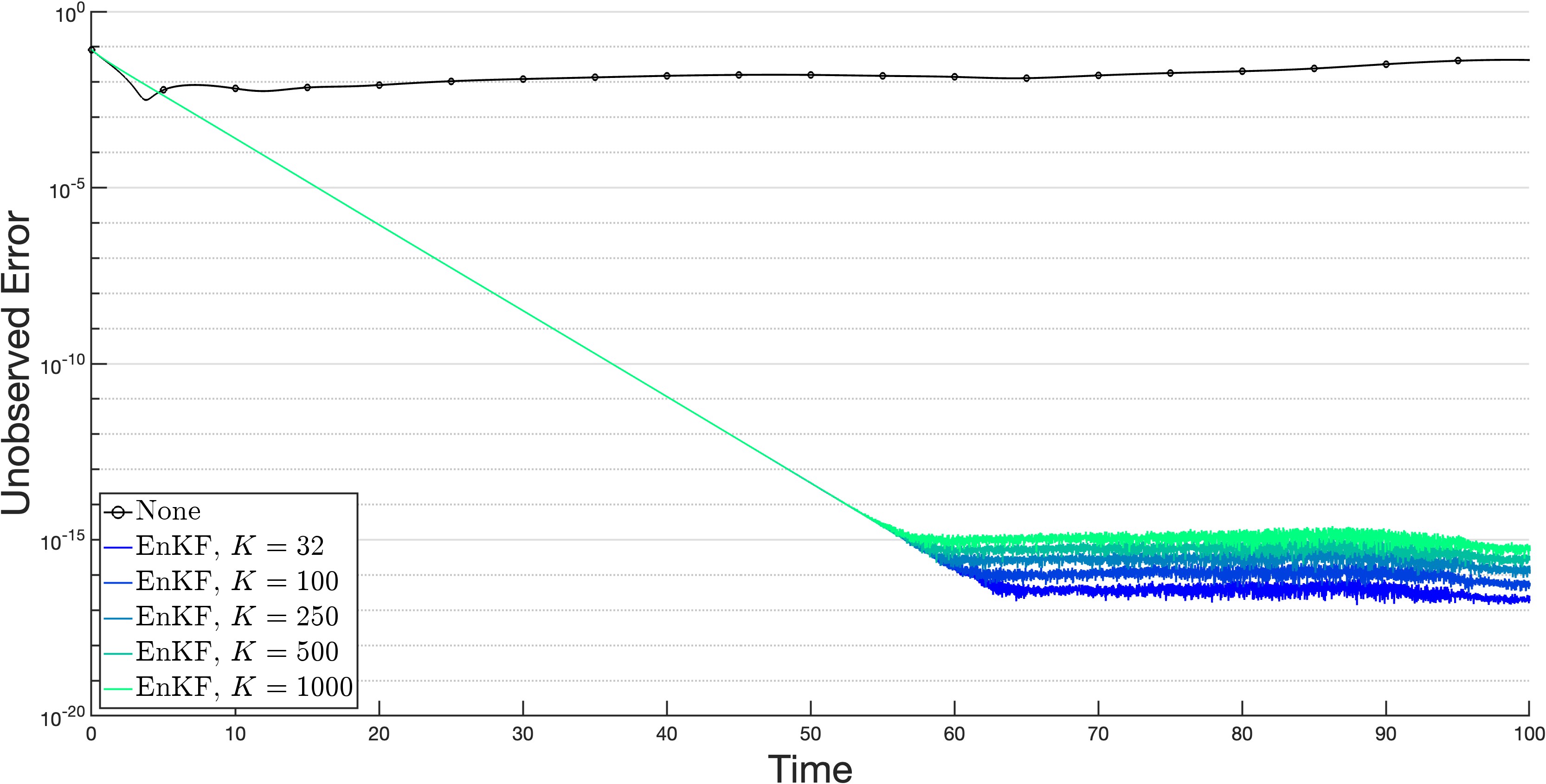}

    \end{subfigure}
    \begin{subfigure}[b]{.32\textwidth}
\centering

            \includegraphics[width=\textwidth]{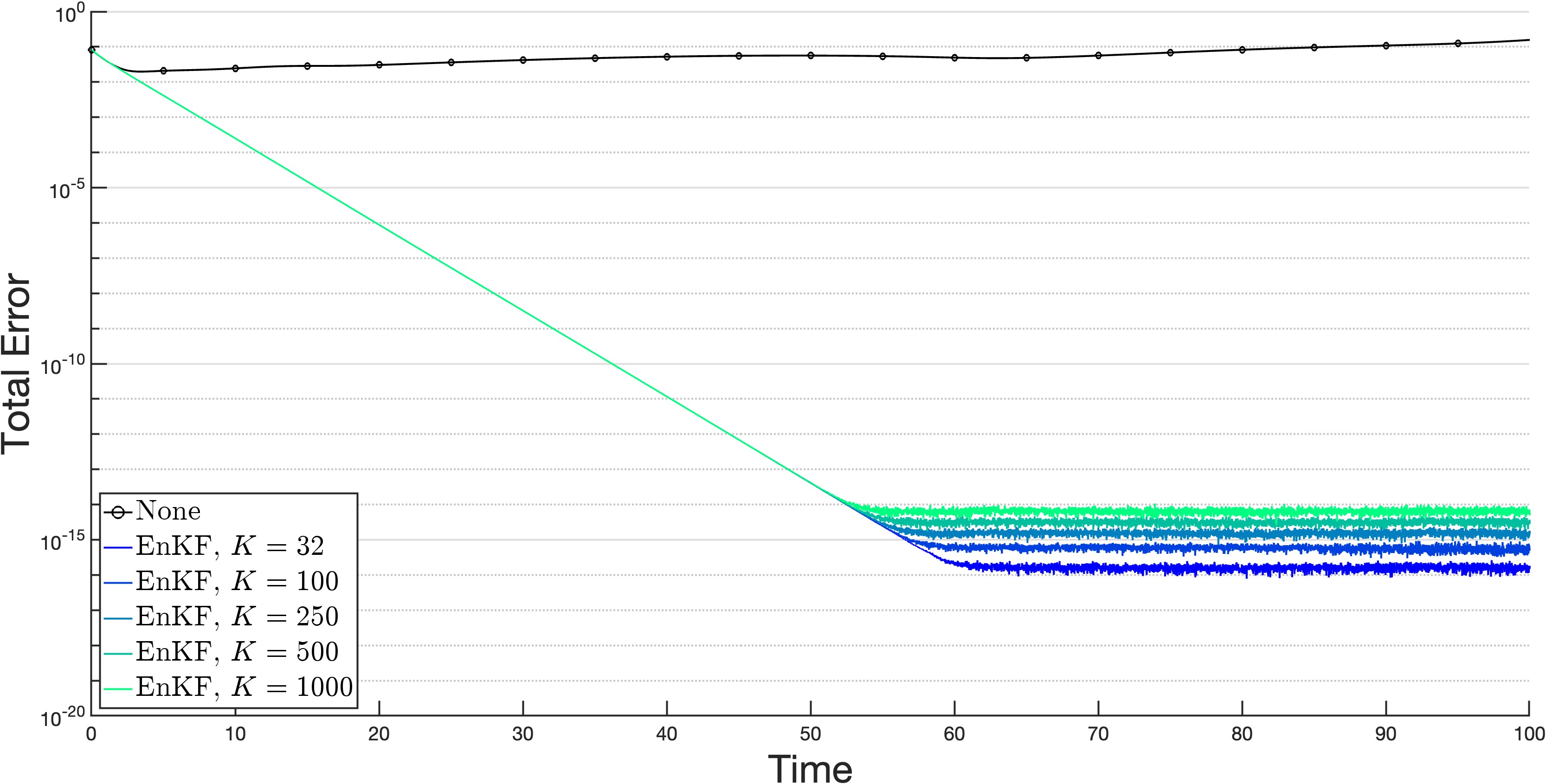}

    \end{subfigure}
    \caption{Error over time for different values of $K$ with covariance inflation. Simulations utilized $\sigma_I^2 = 1e-14$, $\sigma_O^2 = 0$, and $\sigma_E^2 = 1e-16$. }\label{fig:Ensembles_inflation}
\end{figure}

We found that the value of $\sigma_I$ plays a crucial role in the behavior of different numbers of ensemble members, $K$. Figures \ref{fig:Ensemble_variance} and \ref{fig:Ensembles_no_inflation} illustrate the markedly different behavior of the error with a small value of 
$\sigma_I$ compared to $\sigma_I = 0$. Notably, for $\sigma_I = 1e-14$ the behavior is similar for all $K \geq 32$, although larger $K$ values tend to converge to slightly higher error levels. This trend is evident in both Figures \ref{fig:Ensembles_no_inflation} and \ref{fig:Ensembles_inflation} in the long-term error levels for $K = 500$ and $K= 1000$. The reason for this behavior is currently unknown, however the error levels do appear to coincide with the analysis in \cite{Biswas_Branicki_2024} for the EnKF. We note that when observing $M = 16$ modes, it is practically necessary to have $K \geq 2 \cdot M = 32$ ensemble members. Simulations with fewer ensemble members were insufficient to achieve any level of convergence, as the Kalman gain matrix tended to become degenerate and singular. It appears that $\sigma_I$ has a stabilizing effect and is vital to the success of the EnKF, as noted in \cite{asch2016DA}.

\begin{figure}
    \begin{subfigure}[b]{.32\textwidth}
\centering

            \includegraphics[width=\textwidth]{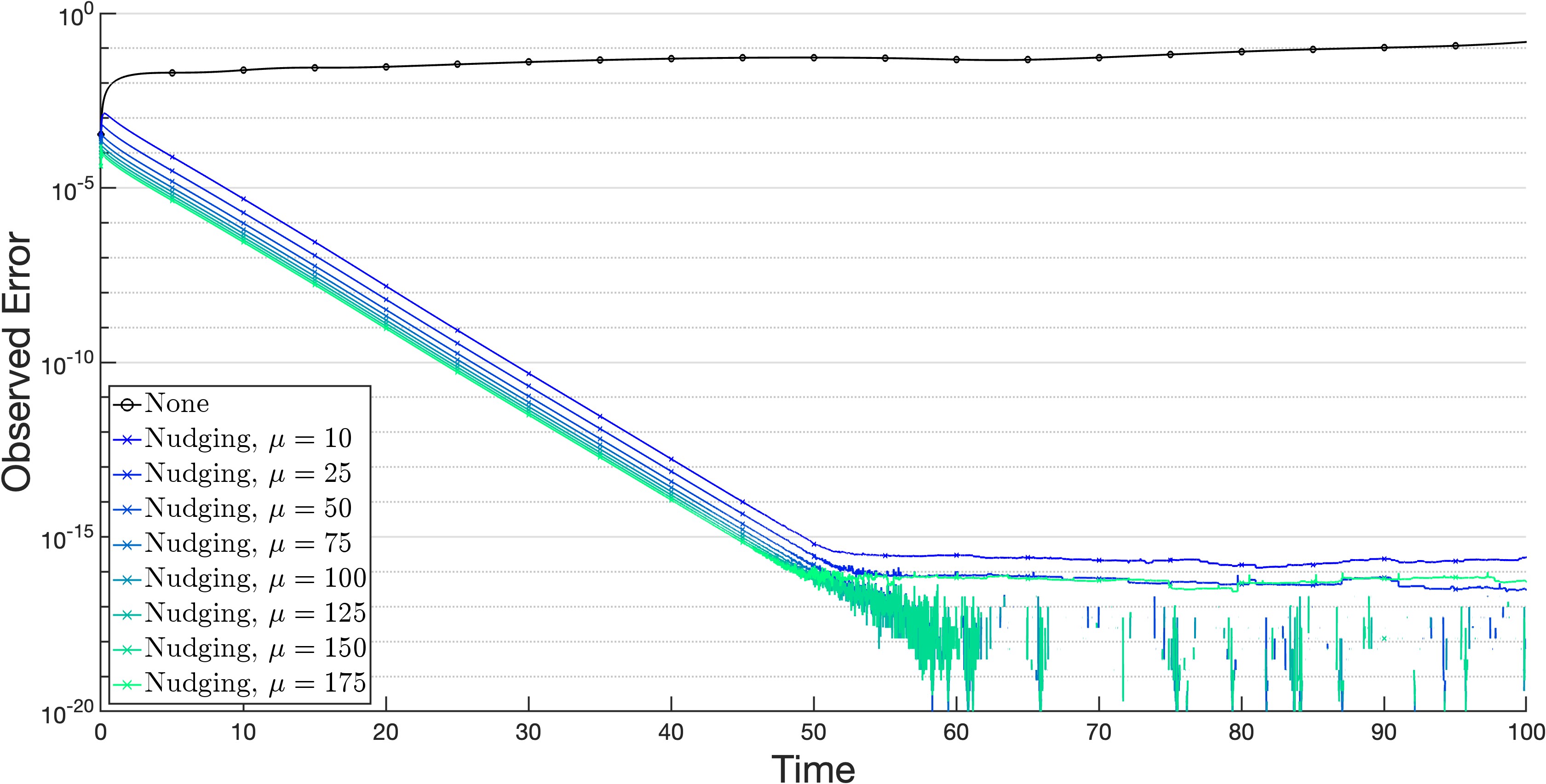}

    \end{subfigure}
    \begin{subfigure}[b]{.32\textwidth}
\centering

            \includegraphics[width=\textwidth]{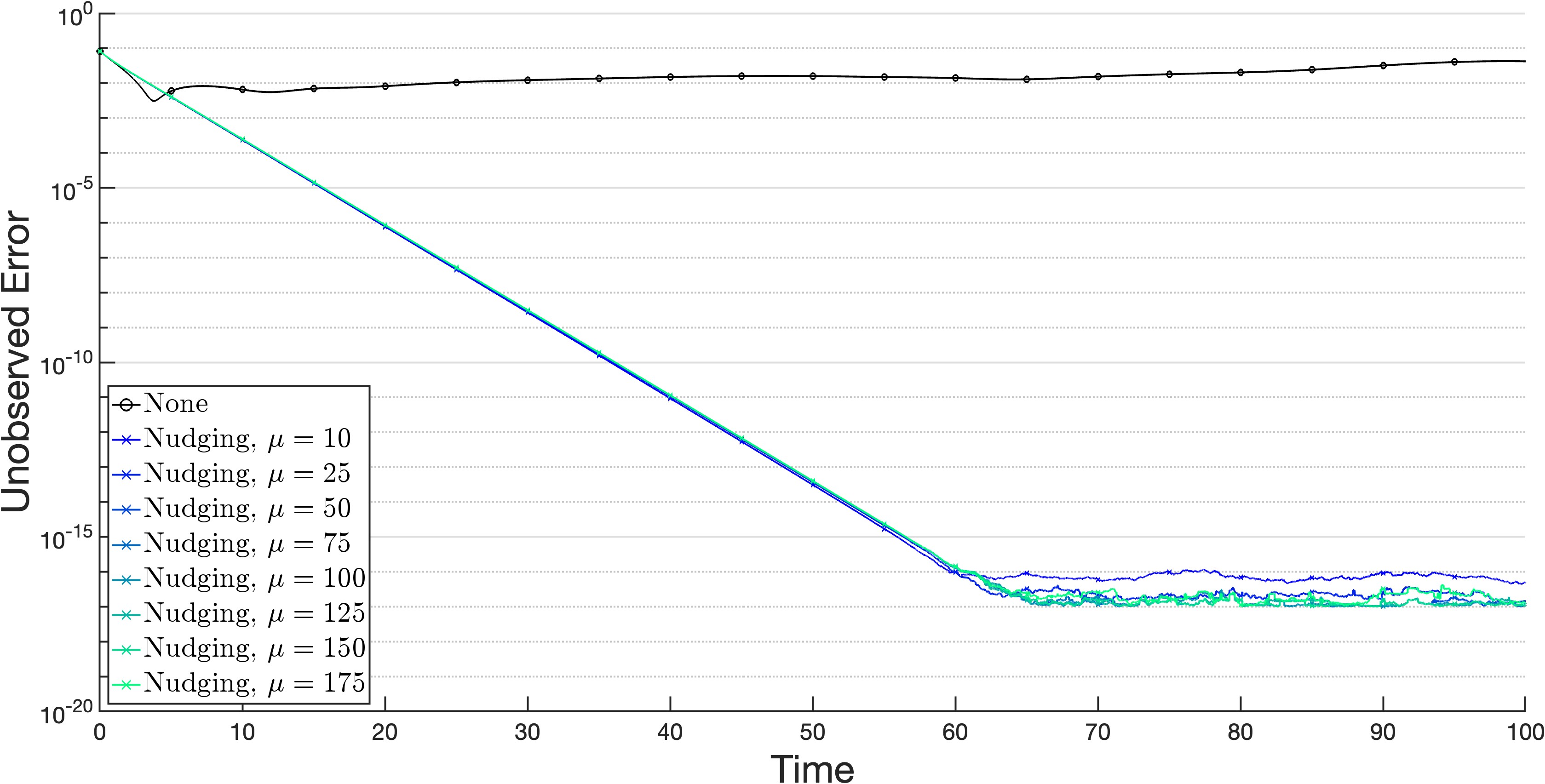}

    \end{subfigure}
    \begin{subfigure}[b]{.32\textwidth}
\centering

            \includegraphics[width=\textwidth]{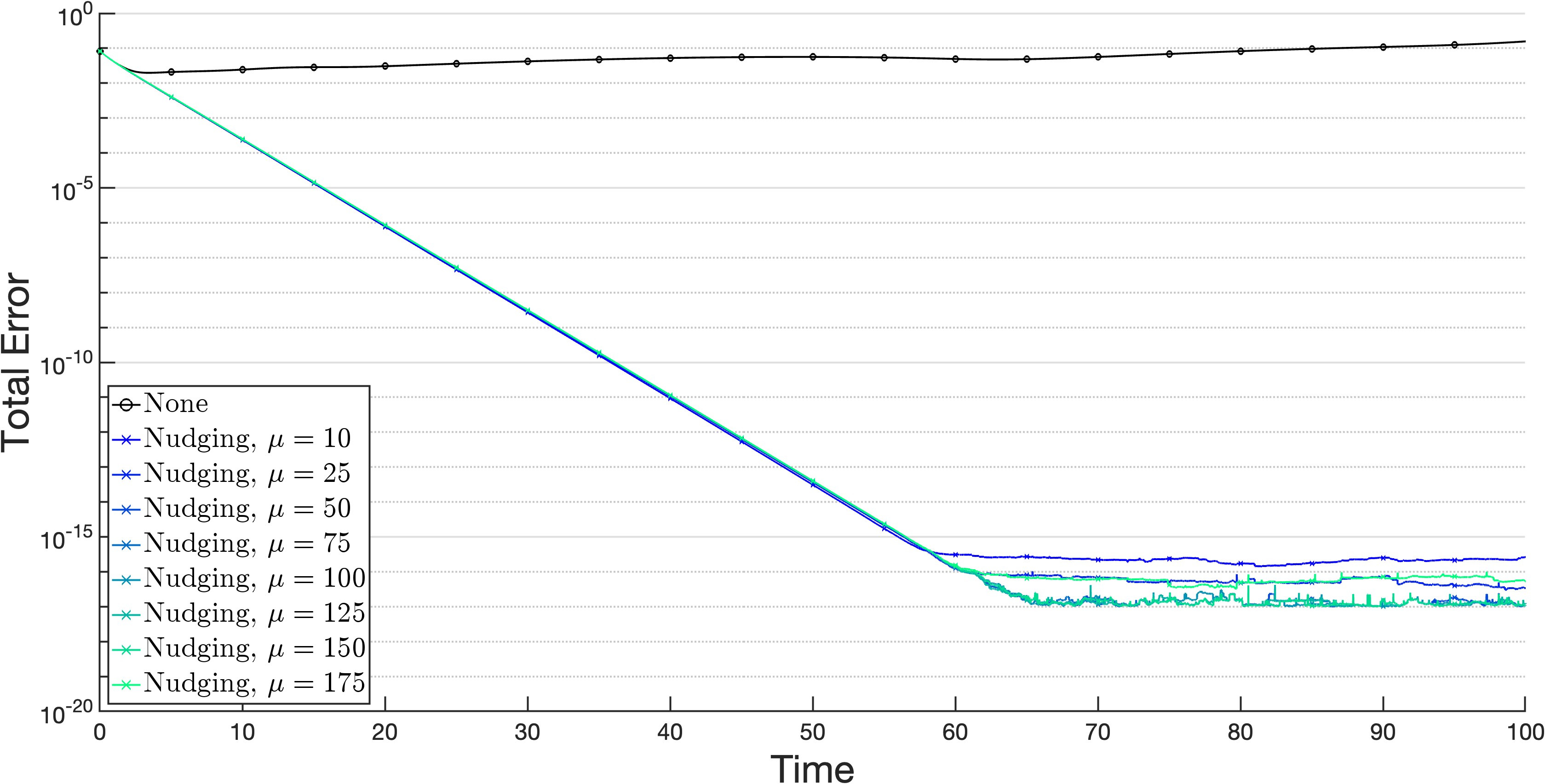}

    \end{subfigure}
    \caption{Error over time for different values of $\mu$ with $\sigma_O^2 = 0$.}\label{fig:mu_no_noise}
\end{figure}

\begin{figure}
    \begin{subfigure}[b]{.32\textwidth}
\centering

            \includegraphics[width=\textwidth]{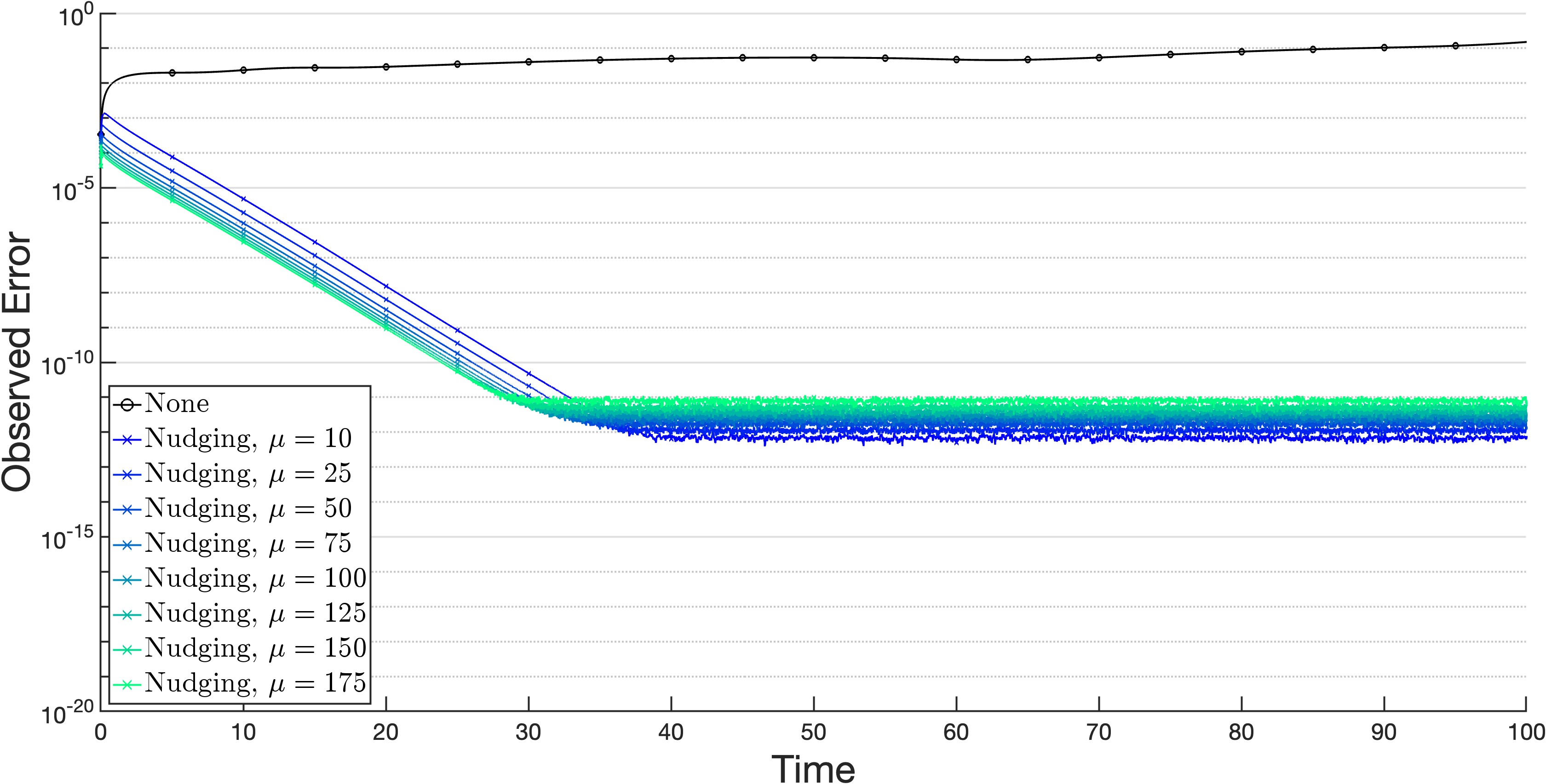}

    \end{subfigure}
    \begin{subfigure}[b]{.32\textwidth}
\centering

            \includegraphics[width=\textwidth]{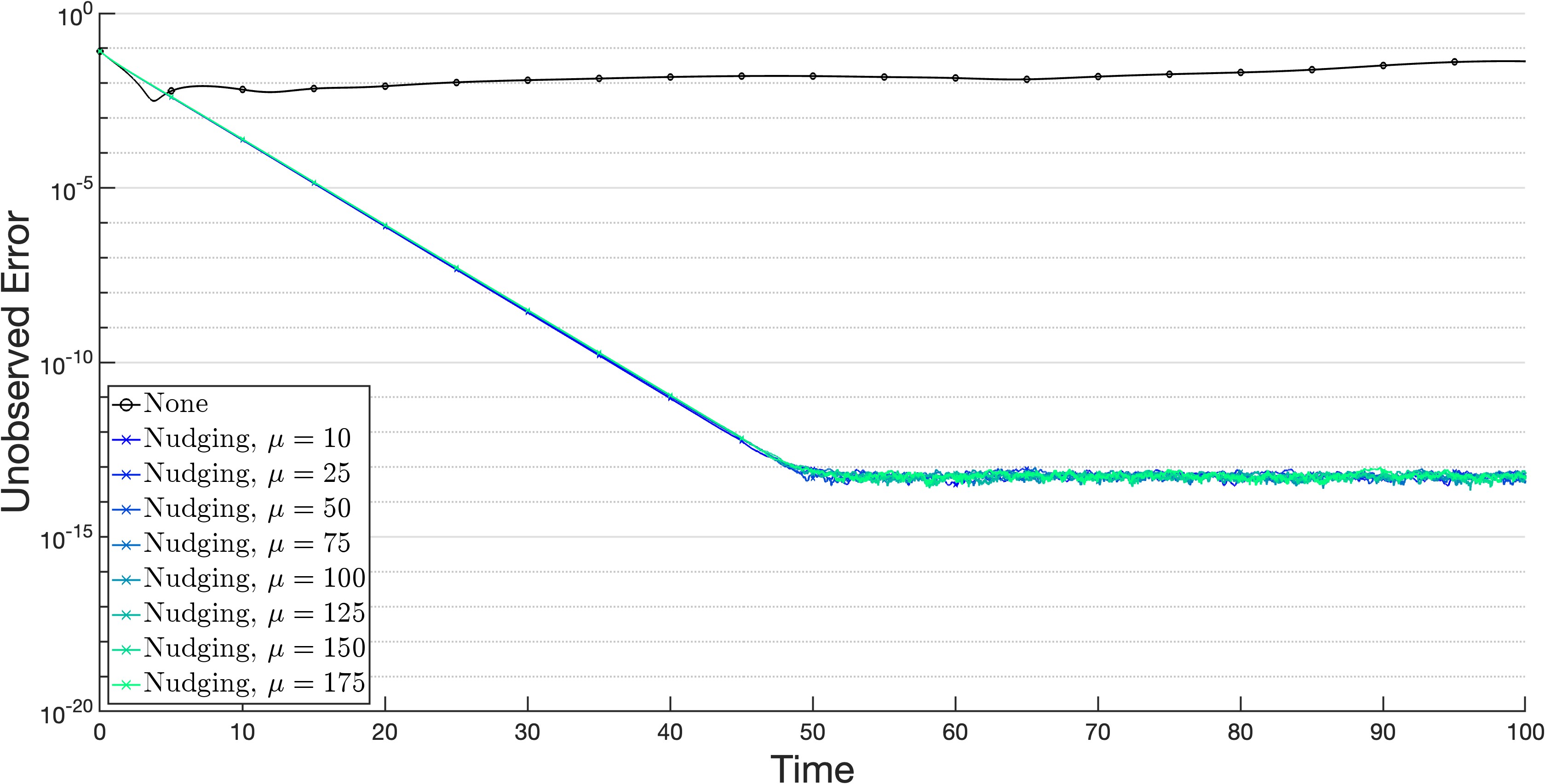}

    \end{subfigure}
    \begin{subfigure}[b]{.32\textwidth}
\centering

            \includegraphics[width=\textwidth]{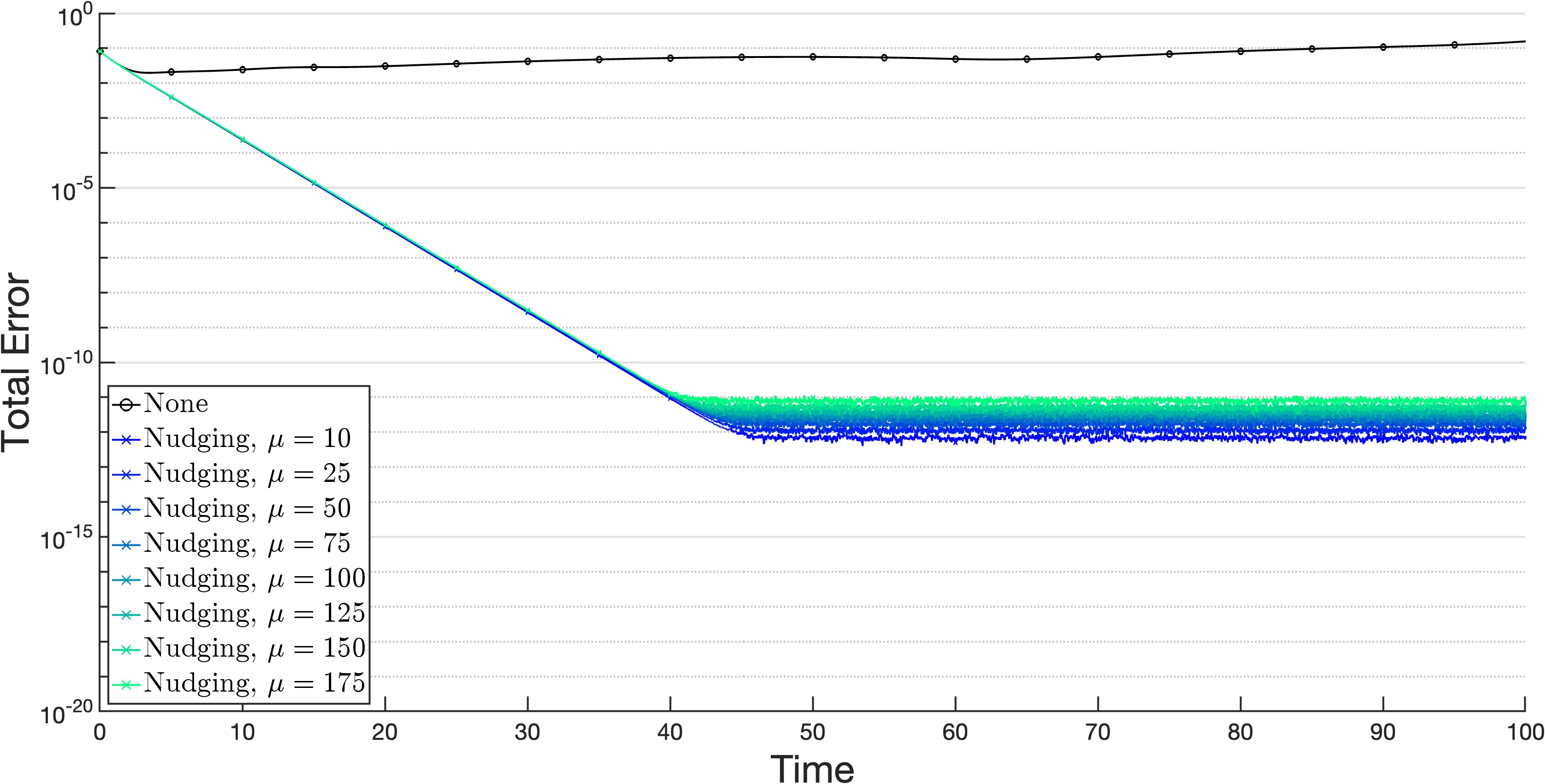}

    \end{subfigure}
    \caption{Error over time for different values of $\mu$ with $\sigma_O^2 = 1e-10$. }\label{fig:mu_noise}
\end{figure}

Lastly, we studied the effect of varying $\mu$ in the presence of noisy observations. Figure \ref{fig:mu_no_noise} shows that the low mode error tends to decrease as $\mu$ increases for $10 \leq \mu \leq 150$; however, the error appears to increase again for $\mu > 150$. This increase is likely due to numerical stability issues, as the feedback-control term is implemented explicitly and thus subject to a CFL condition of 
$\mu \lesssim \frac{2}{\Delta t} = 200$. 
As the value of $\mu$ increases, the numerical scheme becomes unstable, leading to a loss of accuracy and, in severe cases, causing the solution for $v$ blows up numerically. We note that this can be fixed in the code by instead utilizing an implicit implementation, which would relax the CFL condition and allow the use of larger values of $\mu$. However, such a modification is unnecessary, particularly as we are more interested in smaller values of $\mu$ due to the effect that it has on the observational error. Figure \ref{fig:mu_noise} demonstrates that increasing $\mu$ amplifies the error in the observations. We note that this amplification of observational error appears to be consistent with the results of \cite{Bessaih_Olson_Titi_2015}, although the analysis there was for the 2D  incompressible NSE. We note that the effect of the $\mu$ values on convergence rate with observational error will be investigated in a forthcoming work \cite{Synchronization_WIP}.

\begin{figure}
    \begin{subfigure}[b]{.32\textwidth}
\centering

            \includegraphics[width=\textwidth]{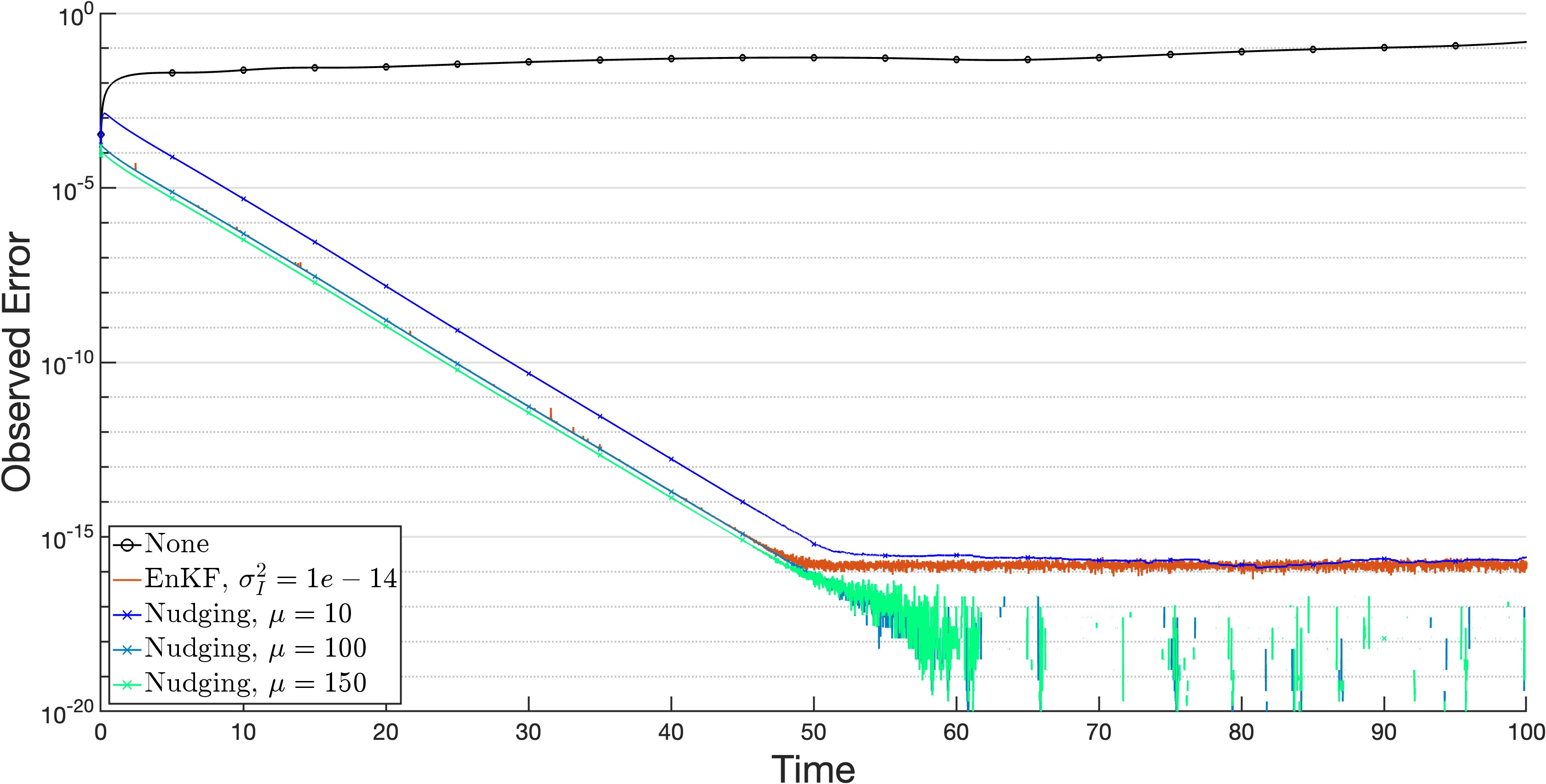}

    \end{subfigure}
    \begin{subfigure}[b]{.32\textwidth}
\centering

            \includegraphics[width=\textwidth]{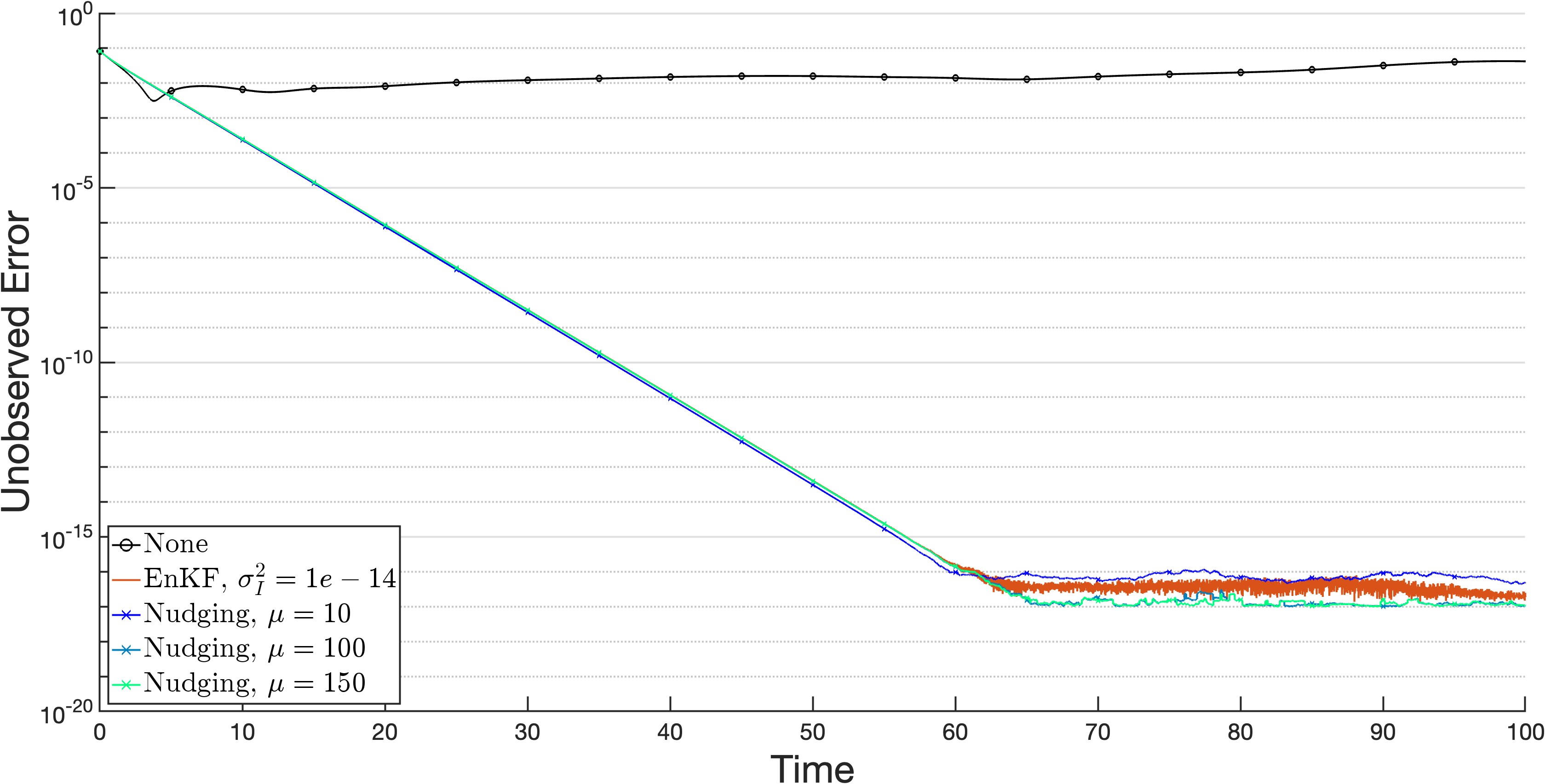}

    \end{subfigure}
    \begin{subfigure}[b]{.32\textwidth}
\centering

            \includegraphics[width=\textwidth]{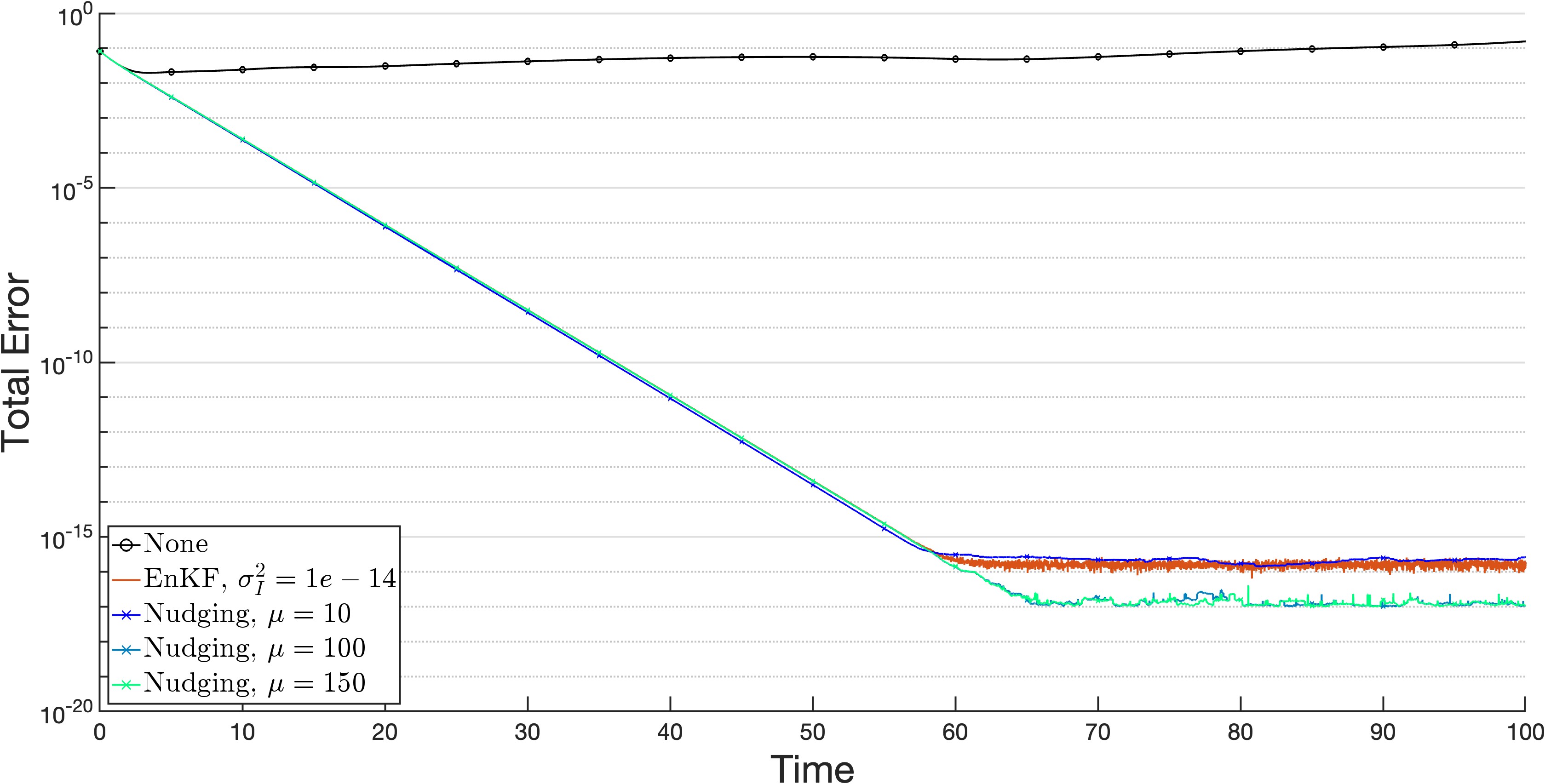}

    \end{subfigure}
    \caption{Error over time comparison of nudging and EnKF.}\label{fig:comparison}
\end{figure}

\begin{figure}
    \begin{subfigure}[b]{.32\textwidth}
\centering

            \includegraphics[width=\textwidth]{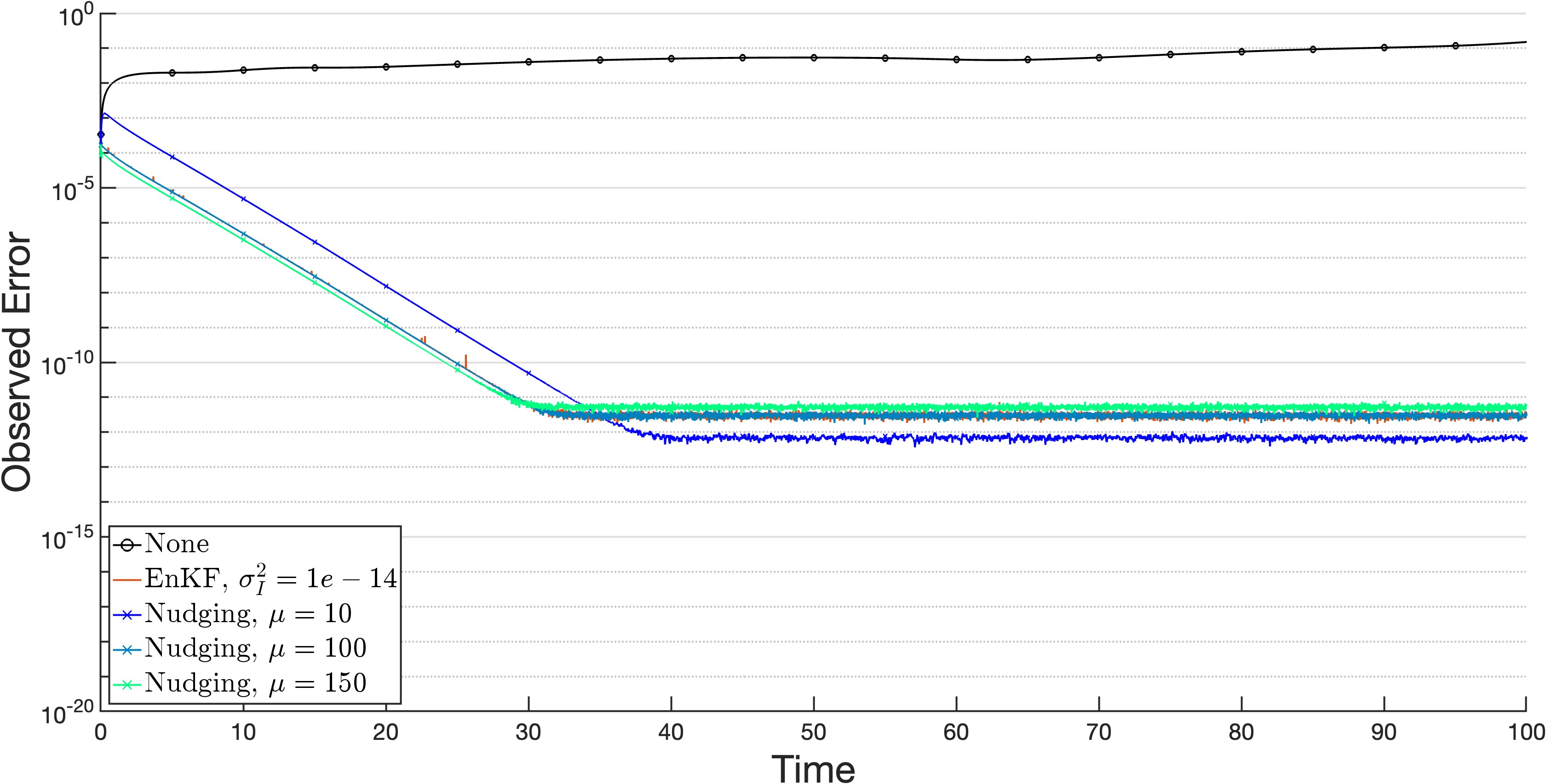}

    \end{subfigure}
    \begin{subfigure}[b]{.32\textwidth}
\centering

            \includegraphics[width=\textwidth]{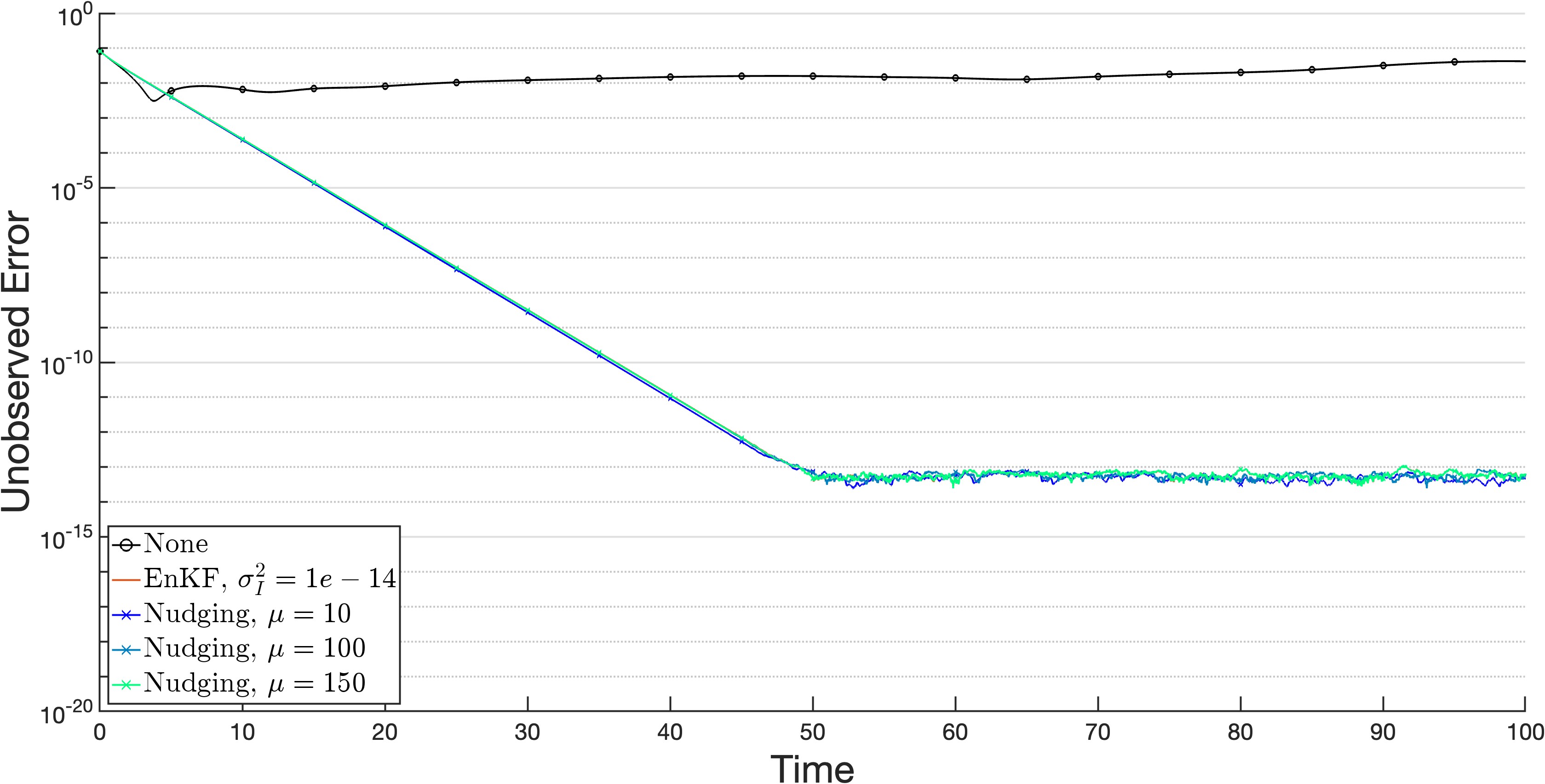}

    \end{subfigure}
    \begin{subfigure}[b]{.32\textwidth}
\centering

            \includegraphics[width=\textwidth]{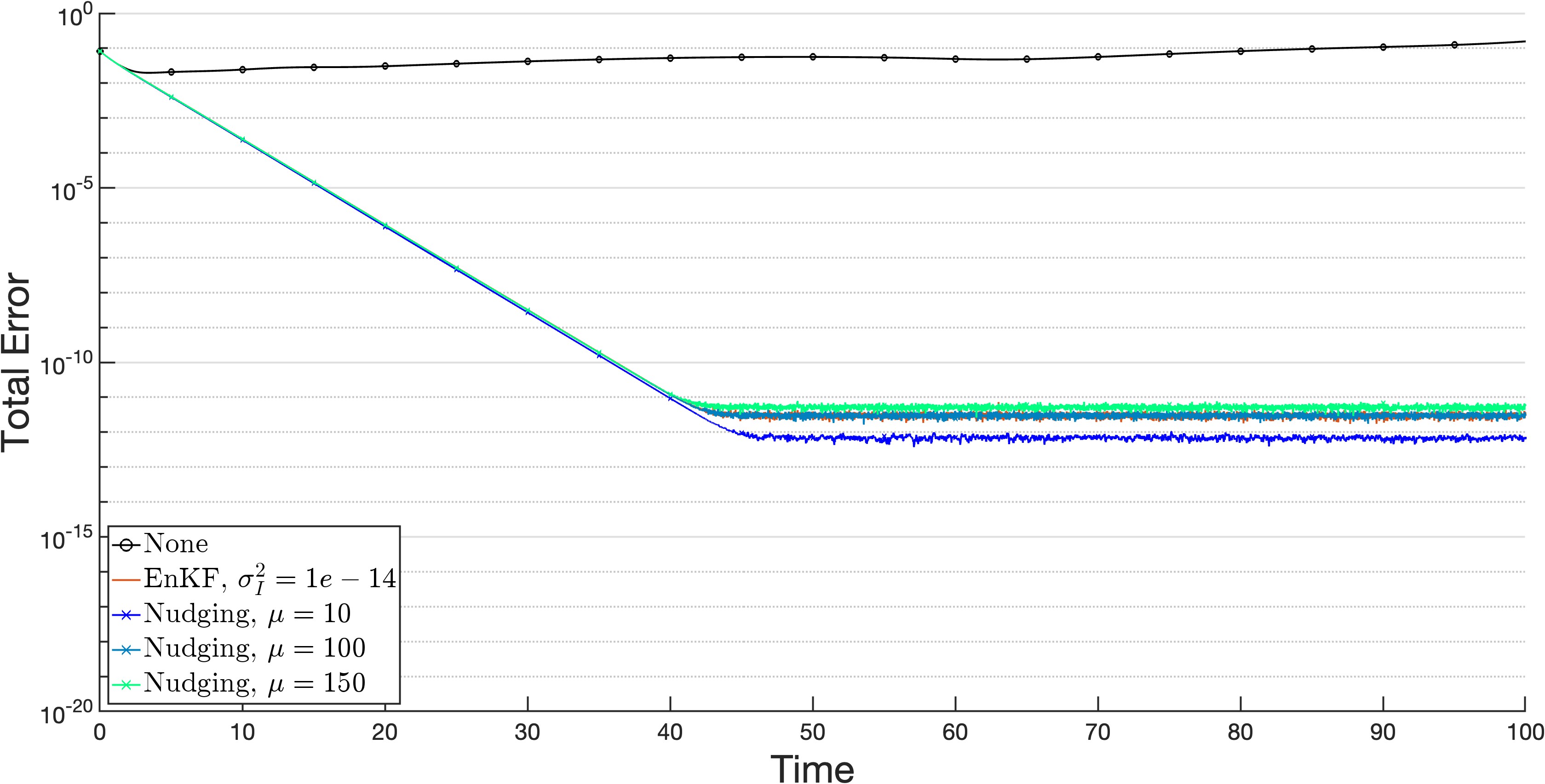}

    \end{subfigure}
    \caption{Error over time comparison of nudging and EnKF.}\label{fig:comparison noisy}
\end{figure}

So far we have discussed the effect of the various choices of parameters, we now compare the EnKF and nudging methods. We have selected the optimal parameters for each method and plotted the results in \Cref{fig:comparison,fig:comparison noisy}. We note that the choice of $\mu$ is crucial for nudging, particularly in the presence of observational noise. Therefore, we present several trials with different $\mu$ values alongside the best choice of parameters for EnKF. The optimal parameters for EnKF appear to be $\sigma_I^2 = 1e-14$, $\sigma_E^2 = 1e-16$, and $K = 32$, although the choice of the number of ensembles is somewhat arbitrary, as indicated by  \Cref{fig:Ensembles_inflation}. In \Cref{fig:comparison}, we observe that both data assimilation methods converge to approximately the same level of error at a similar exponential rate with $\mu = 100$. A similar pattern is observed in  \Cref{fig:comparison noisy} in the presence of observational noise, except that the specific value of $\mu$ affects the inflation of noise in the observations. In both cases, using $\mu = 100$ with nudging achieves almost the exact performance as the EnKF, while using $\mu = 10$ results in slightly slower convergence in the observed error but attains a better overall precision by the end of the simulations.

\label{sec:1d_comparison_results}

\section{Two-dimensional Incompressible Navier-Stokes Equation}
\label{sec:2D}
We now turn our attention to the 2D incompressible NSE, given on the periodic torus $\Omega = \nT^2$, expressed as
\begin{align}\label{NSE}
    u_t - \nu \lap u + u\cdot \nabla u  + \nabla p = f \text{ in } \Omega,\\
    \div{u} = 0\quad\text{and}\quad
    u(0,x) = u_0.\nonumber
\end{align}
Here, $u(t,x)$ denotes the velocity, $f$ is a prescribed body force, $p$ represents the pressure, and $\nu$ is the kinematic viscosity. Note that the body force $f$ is assumed to be constant in time. The corresponding nudging equation for this system is given as follows:
\begin{align}\label{NSE-aot}
    v_t - \mu \lap v + v\cdot \nabla v + \nabla q = f + \mu (P_M \tilde{u} - P_M v) \text{ in } \Omega,\\
    \div{v} = 0\quad\text{and}\quad
    v(0,x) = v_0.\nonumber
\end{align}
Note the three key differences: $q$ represents the pressure for the system of $v$, there is the addition of the feedback-control term $\mu(P_M \tilde{u} - P_M v)$, and the initial condition is given by $v_0$ which is generally arbitrary but primarily taken as 
$v_0 \equiv P_\lambda P_M(\tilde{u}_0)$, where $P_\lambda$ is the Leray projector, which projects onto divergence-free vector fields. 

\subsection{Numerical Methods}

The simulations for this study were conducted on the 2D NSE in the vorticity-stream function formulation using a fully dealiased pseudospectral code, with the physical domain $\Omega = [-\pi, \pi]^2$. Specifically, the 2D NSE \eqref{NSE} is written in the following formulation:
\begin{align}
\psi_t + \triangle^{-1}(u\cdot\nabla\omega) &= \nu\triangle\psi + \triangle^{-1}\nabla^{\perp}\cdot f,\\
u := \nabla^{\perp}\psi = \binom{\partial_y\psi}{-\partial_x\psi}
\quad\text{and}\quad
\omega &:= \partial_{x}u_2-\partial_yu_1 = -\triangle \psi.
\end{align}
Here, the inverse Laplacian $\lap^{-1}$ is taken with respect to periodic boundary conditions and the mean-free condition. For an overview of pseudospectral methods, see, e.g., \cite{Canuto_Hussaini_Quarteroni_Zang_2006,Peyret_2013_spectral_book,Shen_Tang_Wang_2011}. 
While we utilized pseudospectral methods for computing the spatial derivatives, we employed the ETDRK4 scheme described in \cite{Cox_Matthews_2002_ETD} to evolve the equations forward in time. We chose this scheme specifically because it was the time integrator used in \cite{Kelly_Law_Stuart_2014}, whose solution $u$ we aimed to replicate as closely as possible. We note that we use the ETDRK4 scheme to evolve the base NSE equations forward in time, but we implemented the nudging feedback-control term explicitly with forward Euler. The details of our implementation are given in \Cref{alg:ETDRK4}, with $L$ given in this context to be the negative Laplacian operator computed in Fourier space.

In our simulations, we utilized a time step of $\Delta t = 0.1$ and a spatial resolution of $N = 2^7$, which, after dealiasing and consolidating positive and negative frequencies, yields approximately 42 active Fourier modes. We set $\nu = 0.01$ and chose the forcing $f$ such that $ \triangle^{-1}\nabla^{\perp}\cdot f =  f_0\cos(\pi k_f \cdot \vec{x})$, with $k_f = \binom{5}{5}$ and $f_0 = 50$. We chose $f_0 = 50$ to ensure that our parameters correspond to a chaotic regime with a Grashof number $G = 500,000$.
The Grashof number is given by \begin{align} G := \frac{\|f\|_{L^2}}{\nu^2\lambda_1}, \end{align} 
where $\lambda_1$ is the first eigenvalue of the Stokes operator $A$. 
The initial data $u_0$ was generated by evolving initial data $u_0 \equiv 0$ out to time $100,000$ until the resulting solution was approximately on the global attractor.

\begin{figure}
    \centering
    \includegraphics[width=0.75\linewidth]{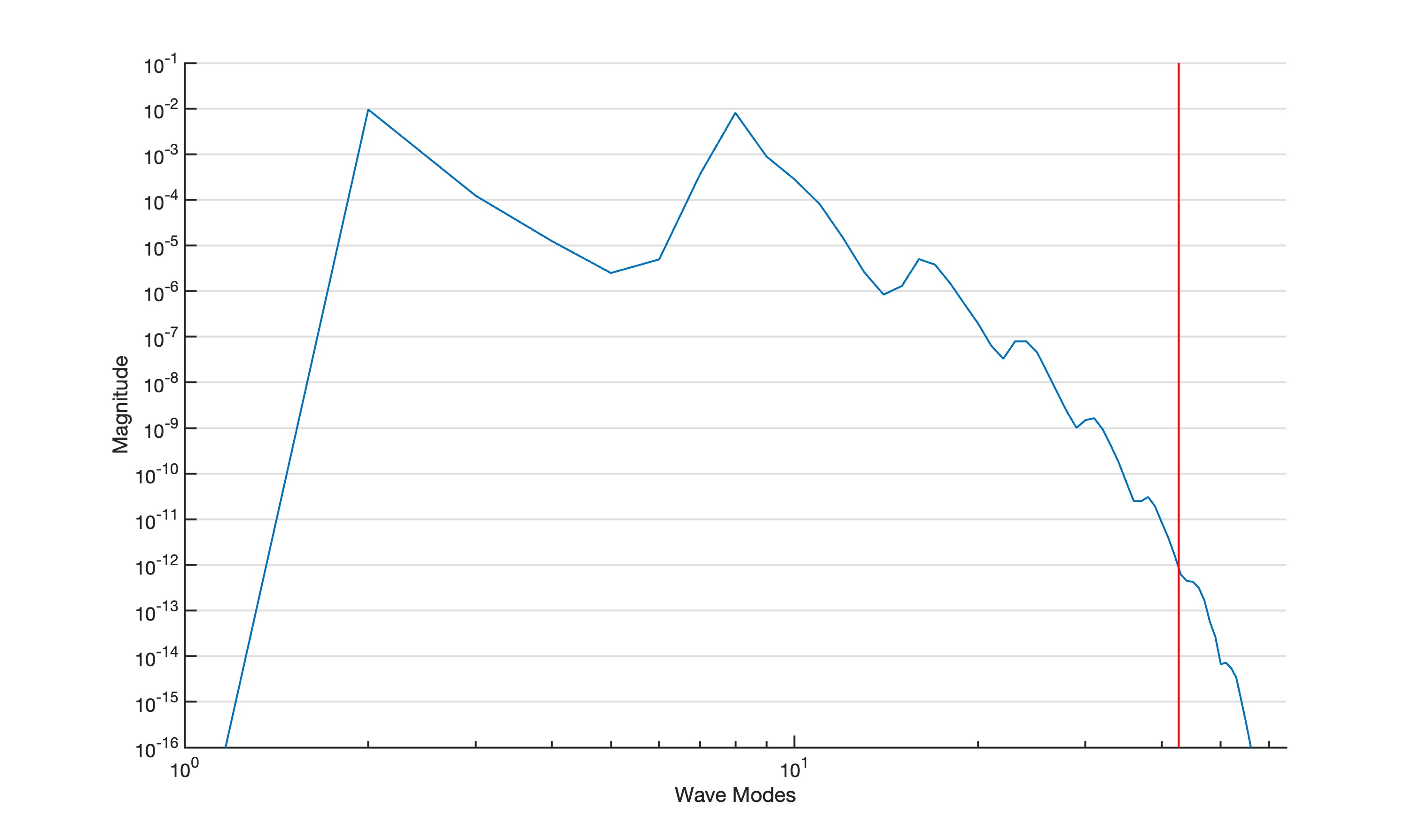}
    \caption{Energy spectrum of the initial data with $\nu = 0.01$, $G = 500,000$, and $\Delta t = 0.01$. The vertical red line is the 2/3 dealiasing cutoff as $\frac{2}{3}\frac{N}{2} = 42.\overline{6}.$}
    \label{fig:spectrum NSE}
\end{figure}

We note that while we utilized the same numerical scheme as that of \cite{Kelly_Law_Stuart_2014}, there are some differences in our approach. Specifically, \cite{Kelly_Law_Stuart_2014} used 3/2 padding for dealiasing, whereas we used 2/3 truncation. Additionally, we did not use the same parameters. In \cite{Kelly_Law_Stuart_2014}, they set 
$\abs{f_0} = 10$, resulting in a Grashof number of 
$G \approx 100,000$, which is one-fifth of our value. Initially, we attempted to replicate their results with lower of these values, but found that the forcing was too weak to induce a chaotic regime. Instead, the system appeared to be in a pseudo-steady state determined primarily by the forcing, with uninteresting dynamics. We can see in \Cref{fig:spectrum NSE} the energy spectrum of the initial condition pictured here, which we note is underresolved as the energy spectrum is approximately 1e-12 at the dealiasing cutoff, but is resolved to single precision. We note that the spectrum apprears to be well behaved for all simulations presented in this section and we did not encounter any numerical instability for our choice of parameters.

The interpolant $P_M$ is given by the projection onto the lowest 10 Fourier modes. We found that this particular choice of modes provides a convenient trade-off between the number of ensembles required and the convergence achieved. In practice, we observed that the number of ensembles needed for the EnKF to work—such that the Kalman gain matrix $K$ tended to be nonsingular and well-behaved—corresponds to the number of observables at a given time step. For the 2D NSE, this is equivalent to the number of points on the integer lattice within a circle of radius $M$, scaling quadratically with the number of observed modes 
$M$. Preliminary tests with a larger number of observations exhibited qualitatively similar behavior, but the simulations were prohibitively expensive for running EnKF over large times.

The nudged solution \eqref{NSE-aot} was simulated using the same numerical scheme as \eqref{NSE}, except that the nudging term was computed explicitly using the observations from the previous time step; see Algorithm \ref{alg:ETDRK4} for full implementation details of the nudging scheme. In our simulations, we chose $\mu = 100$ to satisfy the CFL constraint $\mu \lesssim \frac{2}{\Delta t}$ inherent in this explicit implementation of the nudging term. It's worth noting that we can circumvent this constraint by treating the nudging term implicitly or by employing a higher-order explicit formulation. While such approaches may be advantageous in deterministic scenarios, the value of $\mu$ magnifies the impact of observational noise, as analytically demonstrated in \cite{Bessaih_Olson_Titi_2015}, resulting in reduced precision for large $\mu$. The specific choice of $\mu$ significantly influences the global convergence of the nudged system to the true solution, a topic under exploration in forthcoming work \cite{Synchronization_WIP}, which has already inspired modifications to this algorithm, such as the nonlinear algorithm utilized in \cite{Larios_Pei_2017_KSE_DA_NL,Carlson_Larios_Titi_2023}.

As before, we note that as we are using spectral methods, we are calculating the spatial derivatives using Fourier transforms. In particular, the linear operator $L$ in Algorithm \ref{alg:ETDRK4} is computed on the Fourier wavemodes and used to calculate $\hat{u}_{n+1}$. That is, in this context, $L$ is given to be the second derivative terms of $F$, which is $L := -\lap$. Here $\hat{u}$ indicates the Fourier transform of $u$, given by: 
\begin{align}\hat{u}(t,\vec{k}) := \int_{-\infty}^\infty u(t,\vec{x}) e^{-2\pi i \vec{k}\cdot \vec{x}}d\vec{x}. \end{align} 
As before, we suppress the wavemode index $\vec{k} \in \mathbf{Z}^2$ and index only by the time in our formulations in \Cref{alg:ETDRK4} with the understanding that the equations for $\hat{u}_{n+1}$ and $\hat{v}_{n+1}$ must be solved at each timestep for each wavemode $-N/2 < \abs{\vec{k}} \leq N/2$.  
To obtain the solution in physical space, $u(t,x)$, we simply take the inverse Fourier transform:
\begin{align} {u}(t,\vec{x})= \int_{-\infty}^\infty \hat{u}(t,\vec{k}) e^{2\pi i \vec{k}\cdot \vec{x}}d\vec{k}. \end{align}

\begin{algorithm}
    \caption{Fourth Order Exponential Time Differencing Runge-Kutta (ETDRK4) Scheme (\cite{Cox_Matthews_2002_ETD}) with Explicit Nudging Term}
    \label{alg:ETDRK4}
\begin{algorithmic}[1] 
\Require Observational data $\tilde{u}_n:=P_M(u(t_n)) + \varepsilon(t_n)$, where $\varepsilon(t) \sim \mathcal{N}(0,\sigma_O)$,\\
Observational operator $H:= P_M$ (projection onto the lowest $M$ Fourier modes), \\
Nonlinear function $N(u,t)$,\\
Linear function $L$ for integrating factor,\\
${u}_n$ and ${v}_n$, the value of solutions $u$ and $v$ at previous iteration,\\
$\hat{u}_n = \mathcal{F}(u_n)$ and $\hat{v}_n = \mathcal{F}(v_n)$ represent the Fourier transforms of $u_n$ and $v_n$, respectively.
\bigskip

\State Forward time integration $\hat{u}_{n+1} := M(\hat{u}_n).$
\begin{align}
a &= e^{L \frac{\Delta t}{2}}\hat{u}_n + L^{-1}\left(e^{L\frac{\Delta t}{2} - 1}\right)N(\hat{u}_n,t_n) \\
b &= e^{L \frac{\Delta t}{2}}\hat{u}_n + L^{-1}\left(e^{L\frac{\Delta t}{2} - 1}\right)N(a,t_n + \frac{\Delta t}{2}) \\
c &= e^{L \frac{\Delta t}{2}}a + L^{-1}\left(e^{L\frac{\Delta t}{2} - 1}\right) \left(2N(b,t_n) - N(\hat{u}_n, t_n) \right)\\
    \hat{u}_{n+1} &=e^{L\Delta t}\hat{u}_n + (\Delta t)^{-2}L^{-3}\Bigg\{ \left[-4-L\Delta t + e^{L\Delta t} \left(4-3L\Delta t + (L\Delta t)^2\right) \right] N(\hat{u}_n, t_n)\\
    &\hspace{3cm}+ 2\left[ 2+L\Delta t + e^{L\Delta t}(-2 + L\Delta t) \right](N(a, t_n + \frac{\Delta t}{2}) + N(b, t_n + \frac{\Delta t}{2})\\
    &\hspace{3cm}+ \left[-4 - 3L\Delta t - (L\Delta t)^2 + e^{L\Delta t}(4-L\Delta t) \right] N(c, t_n + \Delta t) \Bigg\}
\end{align}

\State Nudging forward time integration
\begin{align}\hat{v}_{n+1} = M(\hat{v}_n) + \mu \Delta t  P_M\left( \tilde{u}_n - \hat{v}_n\right)\end{align}

\end{algorithmic}
\end{algorithm}

\subsection{Results}
\begin{figure}[htb!]
\centering
	\includegraphics[width = .8\textwidth]{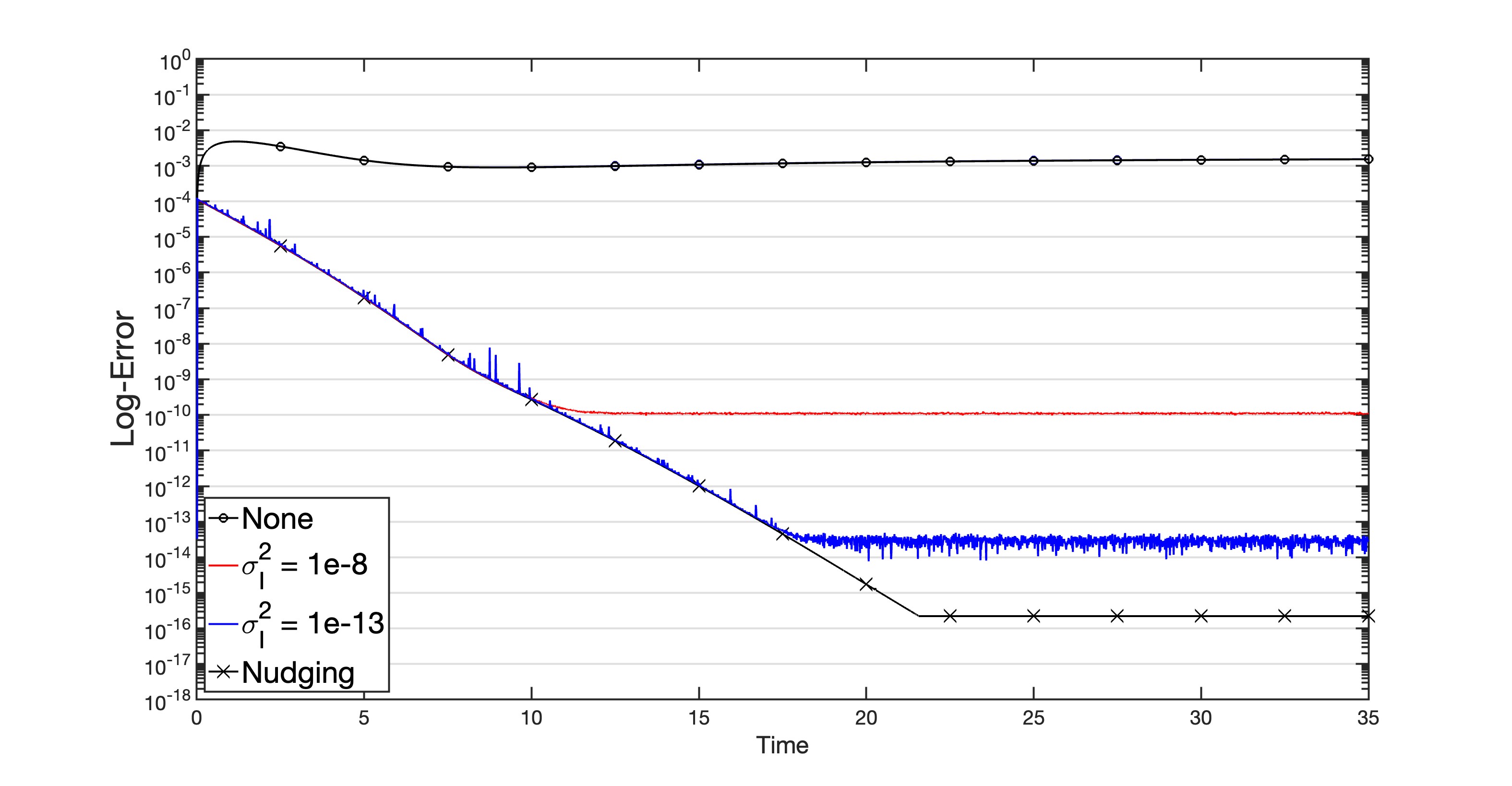}
    \caption{Log-linear plot of $L^2$ error computed on observed modes over time for nudging and EnKF. Lowest 10 Fourier modes were observed with no observational noise. 
}
 \label{fig:low mode determinstic error}
\end{figure}

\begin{figure}[htb!]
\centering
	\includegraphics[width = .8\textwidth]{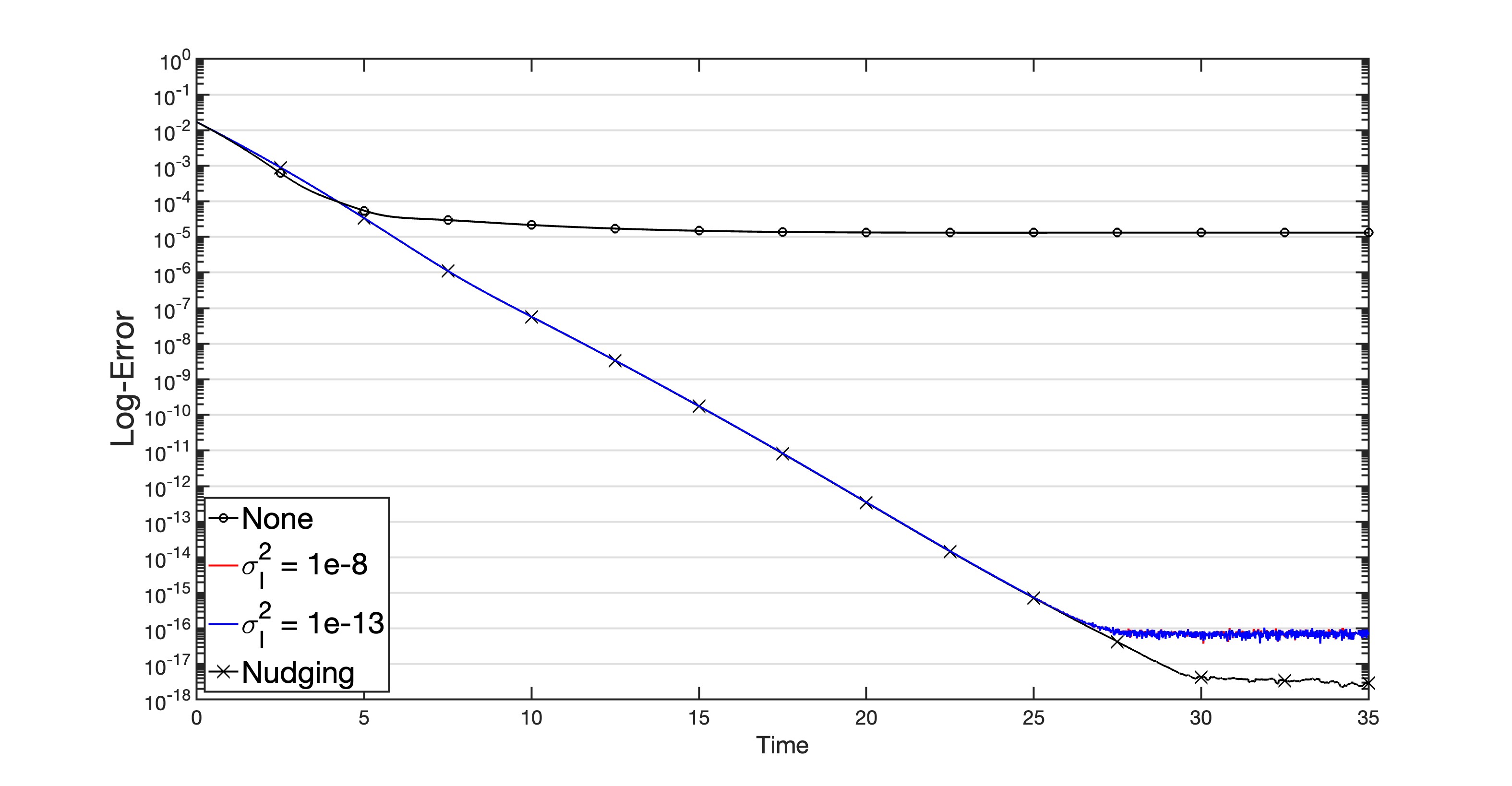}
    \caption{Log-linear plot of $L^2$ error computed on unobserved modes over time for nudging and EnKF. The lowest 10 Fourier modes were observed with no observational noise.}
 \label{fig:high mode determinstic error}
\end{figure}

\begin{figure}[htb!]
\centering
	\includegraphics[width = .8\textwidth]{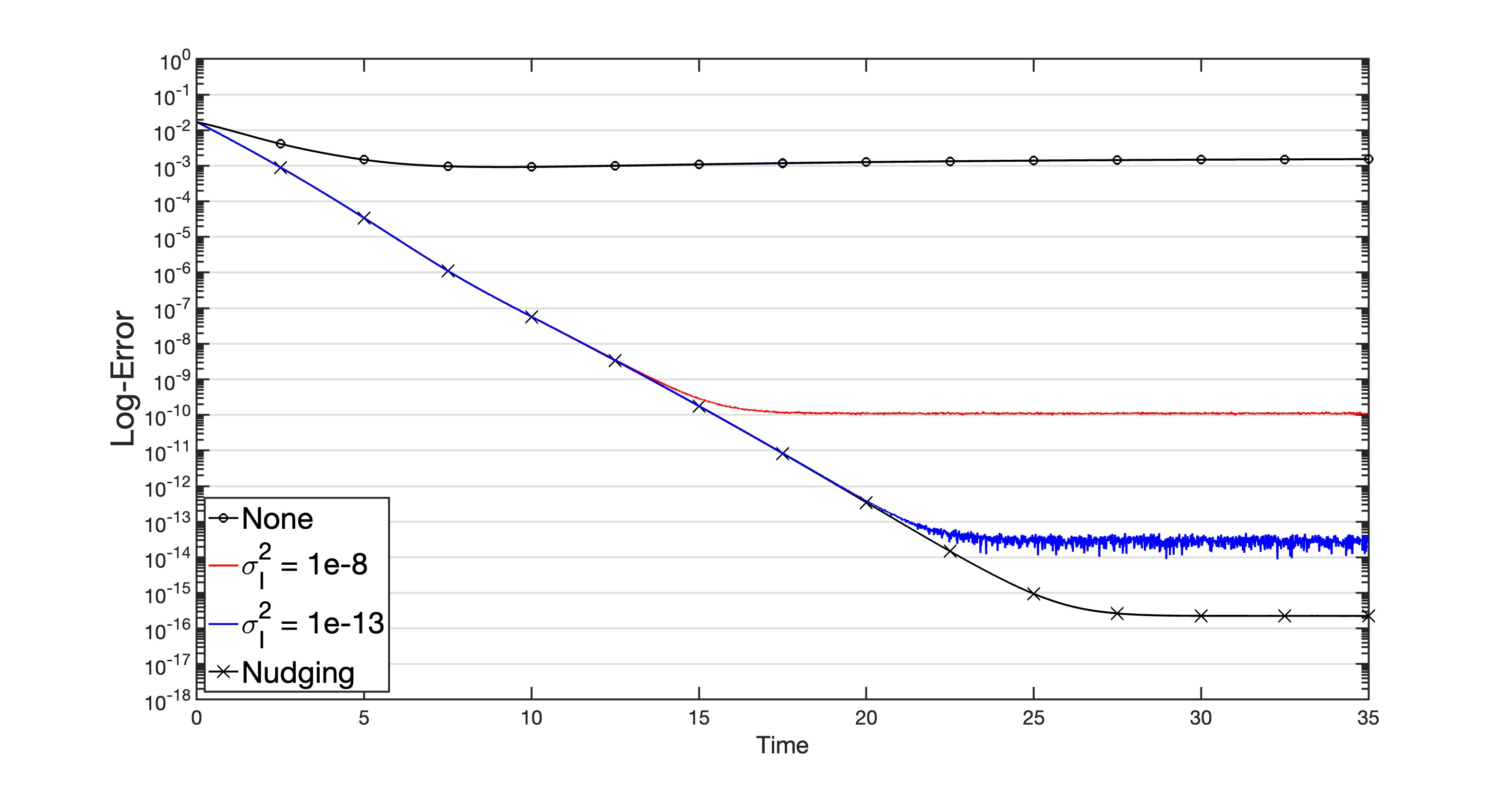}
    \caption{Log-linear plot of $L^2$ error computed on all modes over time for nudging and EnKF. The lowest 10 Fourier modes were observed with no observational noise.}
 \label{fig:all mode determinstic error}
\end{figure}

 To examine these data assimilation algorithms, we used an ``identical twin" experimental design, as illustrated in Section \ref{sec:1d_results}. To test the effectiveness of the EnKF, we first conducted trials with deterministic observations and a fixed ensemble variance $\sigma_E ^2= 1 -14$, while varying the covariance inflation parameter, $\sigma_I$. The results are shown in Figures \ref{fig:low mode determinstic error}-\ref{fig:all mode determinstic error}. The error plotted is the $L^2$ difference between the ensemble average of the solutions produced by the EnKF and the true solution $u$. From Figure \ref{fig:low mode determinstic error}, we observe that all data assimilation methods converge at approximately the same rate to different precision levels. While the convergence rate appears consistent across all methods, the rate of nudging can be influenced by the choice of $\mu$. In contrast, for EnKF, the parameters $\sigma_I$ and $\sigma_E$ seem to affect the level of precision obtained rather than the convergence rate. It is also worth noting that the choice of $\sigma_I$ impacts deviations from the global decay rate. Specifically, a lower value of $\sigma_I$ can sporadically cause large error growth that quickly decays, as evidenced by the spikes on the purple line in Figure \ref{fig:low mode determinstic error}. 

We now turn to examine the error concentrated in the high (unobserved) modes. Figure \ref{fig:high mode determinstic error} shows that both methods of data assimilation converge to machine precision ($2.2204e-16$) at the same rate, which is somewhat surprising. Additionally, the value of $\sigma_I$ does not appear to introduce any artifacts affecting the convergence rate in the high modes. Although nudging seems to achieve a slightly better level of precision, all methods converge within machine precision of the true solution. However, this high-level convergence does not translate to global convergence due to the error inflation in the low modes, as shown in Figure \ref{fig:all mode determinstic error}. It is important to note that $\sigma_I$ is an artificial parameter determined by the user and can be adjusted based on the error in the observed modes. Moreover, one can combine low-mode observations (in the case of deterministic observations) with high-mode approximations produced by the ensembles to eliminate precision loss in the low modes caused by covariance inflation. Thus, in some sense, Figure \ref{fig:high mode determinstic error} represents the true error for the EnKF.


\begin{figure}[htb!]
\centering
	\includegraphics[width = .8\textwidth]{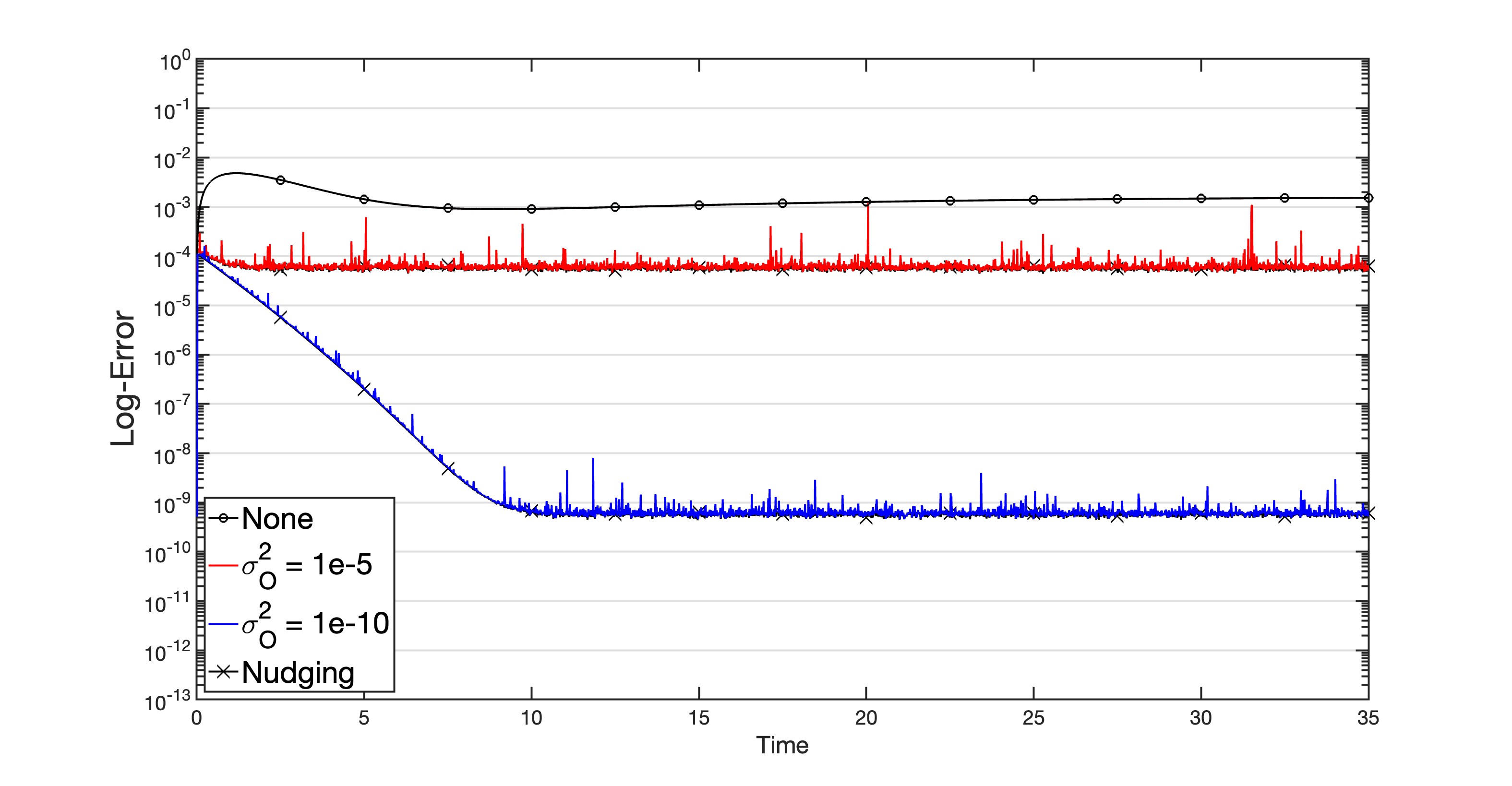}
    \caption{Log-linear plot of $L^2$ error computed on observed modes over time for nudging and EnKF. Coloring, as indicated in legend, corresponds to value of observational noise, $\sigma_O$. Unmarked lines indicate EnKF. In all trials the $10$ lowest Fourier modes were observed.}
 \label{fig:low mode noisy error}
\end{figure}

\begin{figure}[htb!]
\centering
	\includegraphics[width = .8\textwidth]{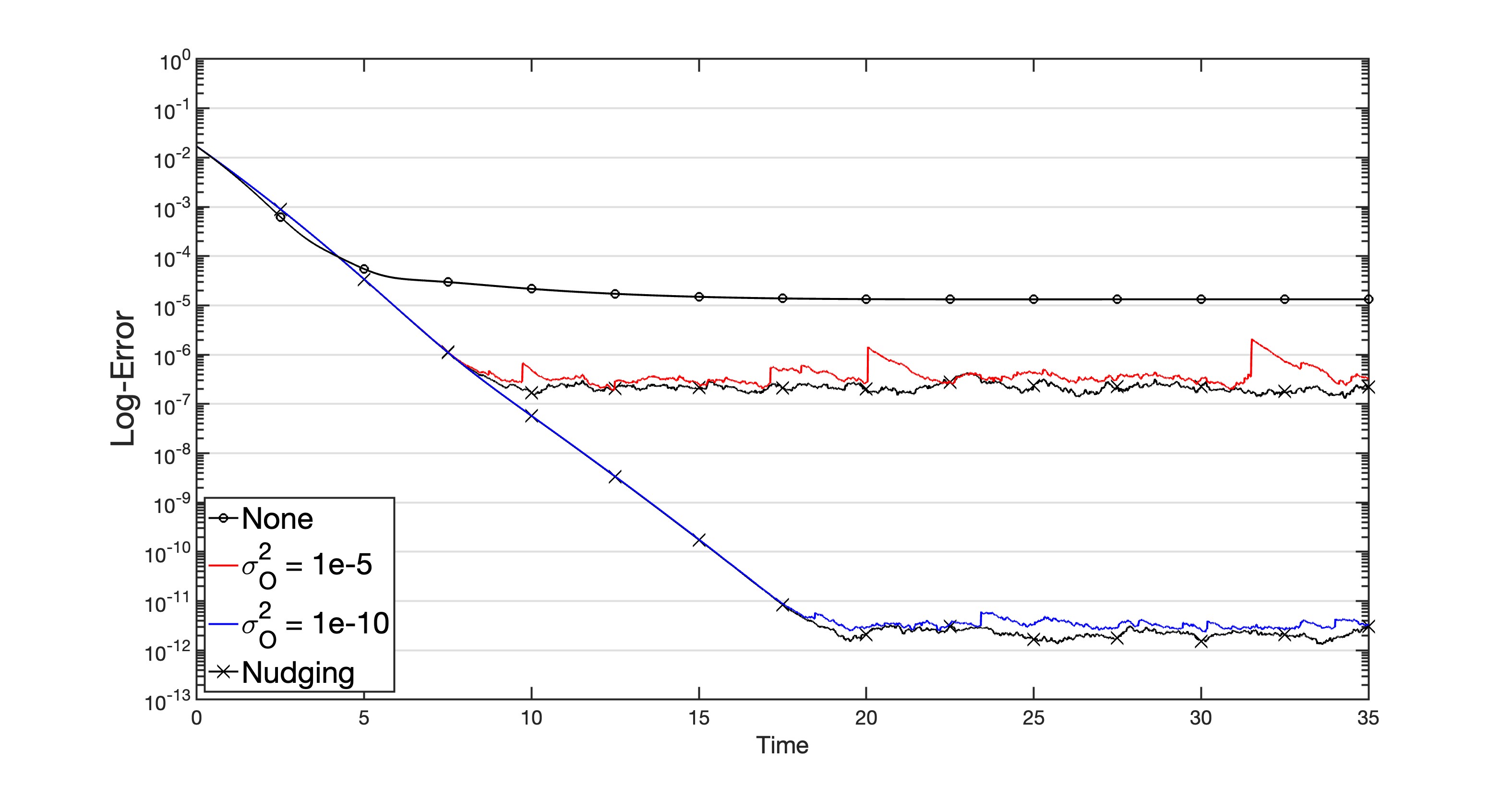}
    \caption{Log-linear plot of $L^2$ error computed on unobserved modes over time for nudging and EnKF. Coloring, as indicated in legend, corresponds to value of observational noise, $\sigma_O$. Unmarked lines indicate EnKF. In all trials the $10$ lowest Fourier modes were observed.}
 \label{fig:high mode noisy error}
\end{figure}

\begin{figure}[htb!]
\centering
	\includegraphics[width = .8\textwidth]{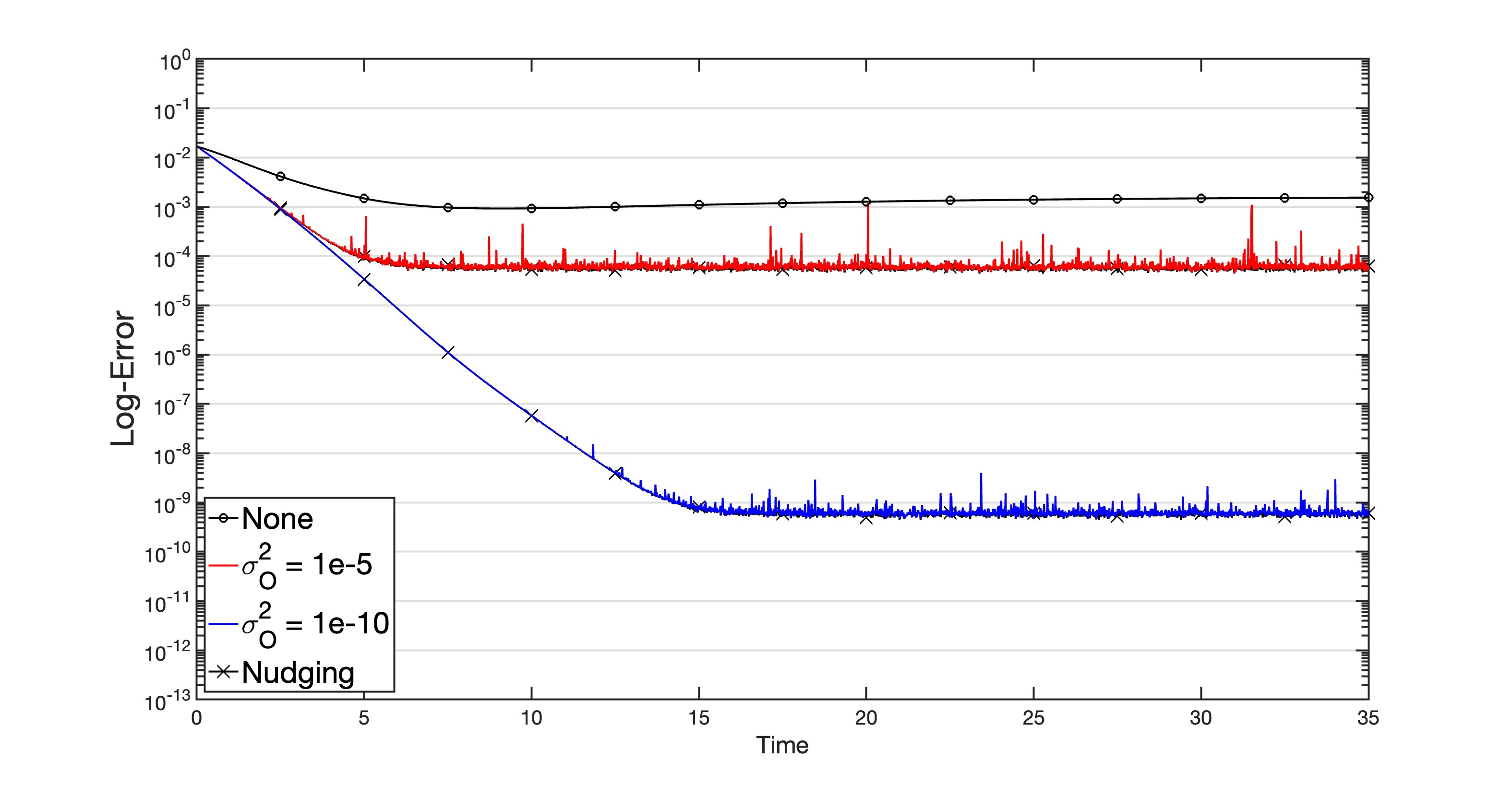}
    \caption{Log-linear plot of $L^2$ error computed on all modes over time for nudging and EnKF. Coloring, as indicated in legend, corresponds to value of observational noise, $\sigma_O$. Unmarked lines indicate EnKF. In all trials the $10$ lowest Fourier modes were observed.}
 \label{fig:all mode noisy error}
\end{figure}

We now move to look at the performance of EnKF and nudging in the presence of observational noise. Holding $\sigma_I$ and $\sigma_E$ constant we look at the performance of EnKF and nudging as we vary the value of $\sigma_O$. We note that we hold $\sigma_I^2 = 1e-13$ and $\sigma_E^2 = 1e-15$ as these appear to give the best levels of convergence for deterministic observations. Figures \ref{fig:low mode noisy error}-\ref{fig:all mode noisy error} reveal several key observations. Firstly, Both nudging and EnKF behave qualitatively the same, converging at essentially the same rate to similar levels of precision (except in the case where $\sigma_I > \sigma_O$). Secondly, The EnKF trials exhibit much more variability in the error, as indicated by the spurious spikes in Figure \ref{fig:low mode noisy error}. In contrast, nudging results in more consistent convergence, with rare deviations from the mean after the error has stabilized to a fixed level of precision. Thirdly, The EnKF experiences transient fluctuations in the error within the high (unobserved) wave modes. The error rapidly increases away from the fixed level of precision and then returns slowly, but still at an exponentially fast rate. This phenomenon is most prominent for $\sigma_O^2 = 1e{-5}$ and corresponds to particularly large deviations in the observed modes. For instance, one such deviation occurs at approximately time $t = 20$, immediately after the error in Figure \ref{fig:low mode noisy error} for $\sigma_O^2 = 1e{-5}$ jumps from $1e{-4}$ to approximately $1 e{-3}$.


\section{Conclusion}
\label{sec:Conclusion}
Through the simulations demonstrated in this work, we have shown computationally that nudging and the EnKF behave similarly under reasonable parameter choices. Specifically, we refer to the parameter $\mu$ for nudging, and $K$, $\sigma_I$, and $\sigma_E$ for the EnKF. Both methods exhibit qualitatively similar behavior, achieving comparable levels of convergence at similar times. While the parameters for each method can be selected to yield similar convergence levels, it is important to note that EnKF is significantly more computationally expensive. This increased cost arises from the need to simulate the governing equations for each ensemble member, scaling with the dimension of the spatial domain. For instance, with only 10 observed modes for the 2D incompressible NSE, EnKF is at least 305 times more expensive than nudging. This computational cost escalates with the number of observations, as we found that with $M$ observed Fourier modes, the required $K$ scales approximately as $M^d$, where $d$ is the spatial dimension of the system. Although parallel simulation of each ensemble member can mitigate this computational burden, assembling the Kalman gain matrix still necessitates global communication between all ensemble members. In general, the primary advantage of EnKF is its ability to reveal statistical information about the system. However, in our computations, we only use the ensemble average to approximate the reference solution, thereby not utilizing this statistical information.

\section*{Acknowledgement}
The research of N.N. was partially supported by NIH grant 1R21AI180492-01 and the Individual Research Grant at Texas A\&M University. 
The research of C.V. was supported in part by NPRP grant \# S-0207-359200290 from the Qatar National Research Fund (a member of Qatar Foundation).

\bibliographystyle{plain}
\bibliography{VictorBiblio}
\end{document}